\documentclass{amsart}

\usepackage{amsfonts,amsmath,amssymb,amsthm}
\usepackage{appendix}
\usepackage{geometry}
\usepackage{graphics}
\usepackage{hyperref}
\usepackage{latexsym}

\theoremstyle{plain}
\def\beq{\begin{eqnarray}}
\def\eeq{\end{eqnarray}}
\def\beqq{\begin{eqnarray*}}
\def\eeqq{\end{eqnarray*}}

\def\trace{\mbox{tr}}

\def\qed{\hfill$\sqcap\kern-8.0pt\hbox{$\sqcup$}$\\}
\def\beq{\begin{eqnarray}}
\def\eeq{\end{eqnarray}}
\def\be{\begin{equation}}
\def\ee{\end{equation}}

\numberwithin{equation}{section}

\newtheorem{theorem}{Theorem}[section]
\newtheorem{assumption}{Assumption}
\newtheorem{claim}[theorem]{Claim}

\newtheorem{corollary}[theorem]{Corollary}
\newtheorem{example}[theorem]{Example}
\newtheorem{lemma}[theorem]{Lemma}
\newtheorem{proposition}[theorem]{Proposition}

\theoremstyle{definition}
\newtheorem{definition}[theorem]{Definition}
\theoremstyle{remark}
\newtheorem{remark}[theorem]{Remark}

\makeatother

\begin{document}

\title{Examples of Einstein manifolds in odd dimensions}

\author{Dezhong Chen}
\address{Department of Mathematics, University of Toronto, Toronto, Ontario, M5S 2E4, Canada}
\email{dchen@math.utoronto.ca}

\subjclass{53C25}
\keywords{Conformally compact Einstein manifold, $Q$-curvature, Kreck-Stolz invariant}

\maketitle

\begin{abstract}
We construct Einstein metrics of non-positive scalar curvature on certain solid torus bundles over a Fano K\"{a}hler-Einstein manifold. We show, among other things, that the negative Einstein metrics are conformally compact, and the Ricci-flat metrics have slower-than-Euclidean volume growth and quadratic curvature decay. Also we construct positive Einstein metrics on certain $3$-sphere bundles over a Fano K\"{a}hler-Einstein manifold. We classify the homeomorphism and diffeomorphism types of the total spaces when the base is the complex projective plane.
\end{abstract}

\section{Introduction}\label{intr}

An \textit{Einstein} manifold is a (pseudo-)Riemannian manifold whose Ricci tensor is proportional to the metric tensor \cite{Bes87}, i.e.,
\begin{align}\label{ee}
\mbox{Ric}(g)=\lambda g.
\end{align}
From the analytic point of view, the Einstein equation \eqref{ee} is a complicated non-linear system of partial differential equations, and it is hard to prove the existence of Einstein metrics on an arbitrary manifold. For example, it is still unknown whether every closed manifold of dimension greater than $4$ carries at least one Einstein metric \cite[0.21]{Bes87}. Thus people turn to studying Einstein manifolds with large isometry group, e.g., when the isometry group acts on the manifold transitively, and hence \eqref{ee} reduces to a system of algebraic equations, or when the isometry group acts on the manifold with principal orbits of codimension one, and hence \eqref{ee} reduces to a system of ordinary differential equations. In both cases, it becomes much more manageable to establish some existence results for Einstein metrics (see the surveys \cite[Chap 7]{Bes87} and \cite[\S\S2,4]{Wan98}). Recent progress in this direction includes the variational approach to study homogenous Einstein metrics by several authors (see \cite{Lau02} in the noncompact case, and \cite{BohWanZil04} in the compact case).

Another natural simplification of \eqref{ee} is to impose the Einstein condition on the total space of a \textit{Riemannian submersion with totally geodesic fibers} \cite[9.61]{Bes87}. One major obstacle in this setup arises from the existence of Yang-Mills connections with curvature form of constant norm. In general, it is still unknown when this necessary condition is satisfied. However, if the structure group of the underlying fiber bundle is a compact \textit{torus} and the base is closed, then the curvature form of a principal connection is the pullback of a harmonic $2$-form on the base. This observation leads people to consider principal torus bundles over products of K\"{a}hler-Einstein manifolds and their associated fiber bundles. In these cases, the harmonic $2$-forms are rational linear combinations of the Ricci forms. Many interesting Einstein metrics have been discovered on these spaces, among which is the first inhomogeneous Einstein metric with positive scalar curvature constructed by Page \cite{Pag78} on $\mathbb{C}P^2\sharp\overline{\mathbb{C}P}^2$, i.e., the nontrivial $2$-sphere bundle over $\mathbb{C}P^1$. Later on, the method of Page was extended by several authors to construct positive Hermitian-Einstein metrics on certain $2$-sphere bundles over products of Fano K\"{a}hler-Einstein manifolds \cite{Ber82, PagPop87, WanWan98}. On the other hand, Koiso and Sakane \cite{KoiSak86} were able to construct the first inhomogeneous Fano K\"{a}hler-Einstein metrics on certain $2$-sphere bundles over two copies of the same Fano K\"{a}hler-Einstein manifolds. Note that all the aforementioned Einstein metrics live on the $2$-sphere bundles associated with a principal \textit{circle} bundle. We refer the reader to \cite[Chap 9]{Bes87} and \cite[\S\S1,3]{Wan98} for more details.

In this paper, we focus on constructing smooth Einstein metrics on the solid torus bundles and the $3$-sphere bundles associated with a principal \textit{$2$-torus} bundle over a single Fano K\"{a}hler-Einstein manifold. Note that our principal $2$-torus bundles are products of $S^1$ with a principal circle bundle (see Proposition \ref{product}).

First we consider the associated solid torus bundles. In general, given a principal $2$-torus bundle $P$, we can construct a solid torus bundle $T(P)=P\times_{S^1\times S^1}(B^2\times S^1)$ via the left action of $2$-torus $S^1\times S^1$ on solid torus $B^2\times S^1$
\begin{align*}
(e^{i\alpha},e^{i\beta})\cdot(re^{ix},e^{iy})=(re^{i(x+\alpha)},e^{i(y+\beta)}).
\end{align*}
Our first result demonstrates the existence of {\it conformally compact} Einstein (CCE) metrics on $T(P)$.

\begin{theorem}\label{neg}
Let $V$ be a Fano K\"{a}hler-Einstein manifold with first Chern class $pa$, where $p\in\mathbb{Z}_+$, and $a\in H^2(V;\mathbb{Z})$ is an indivisible class. For $q=(q_1,q_2)\in\mathbb{Z}\oplus\mathbb{Z}$ with $|q_2|>0$, let $P_q$ be the principal $2$-torus bundle over $V$ with characteristic classes $(q_1a,q_2a)$. Then there exist a two-parameter family of CCE metrics on $T(P_q)$.
\end{theorem}

\begin{remark}\label{pw}
In the case of $q_2=0$ and $|q_1|>p$, we constructed in \cite[Theorem 1.4]{Che09} a one-parameter family of CCE metrics on $T(P_q)$.
\end{remark}

\begin{remark}
All the CCE manifolds Theorem \ref{neg} yields are of odd dimensions, among which the lowest-dimensional ones are some nontrivial solid torus bundles over $\mathbb{C}P^1$.
\end{remark}

\begin{remark}
By Theorem \ref{neg}, the moduli space of CCE structures on $T(P_q)$ is nonempty. Then it must be a smooth infinite-dimensional Banach manifold (see \cite[\S2]{And08}).
\end{remark}

A CCE manifold is a complete noncompact Einstein manifold which is conformal to the interior of a compact Riemannian manifold-with-boundary \cite{Gra00, And05, Che09}. The conformal boundary of the bulk space is called the {\it conformal infinity}. The fundamental link between the global geometry of a CCE manifold and the conformal geometry of its conformal infinity lies in the asymptotic expansion of its volume function w.r.t. the {\it geodesic} defining function determined by a choice of representative for the conformal infinity. In odd dimensions, the expansion has a logarithmic term whose coefficient turns out to be independent of the choice of representative, and is identified with a constant multiple of the total $Q$-curvature of the conformal infinity \cite{GraZwo03, FefGra02}. Moreover, the pointwise $Q$-curvature can be read off from the asymptotic behavior of a certain formal solution to a Poisson equation on the CCE manifold \cite[Theorem 4.1]{FefGra02}. Applying these general results to the CCE manifolds in Theorem \ref{neg} reveals that

\begin{theorem}\label{zq}
The conformal infinity of every CCE manifold in Theorem \ref{neg} has a $Q$-flat representative, i.e., a metric with constant zero $Q$-curvature.
\end{theorem}

\begin{remark}
It is worth mentioning that the same feature holds true for all the odd-dimensional CCE manifolds in \cite{Che09} (cf. \cite[Theorem 1.10]{Che09}).
\end{remark}

\begin{remark}
As an analogue to the well-known Yamabe problem, one may wonder whether every conformal class of metrics on an even-dimensional closed manifold has a representative with constant $Q$-curvature. We refer the reader to \cite{ChaYan95, Bre03, Ndi07, DjaMal08} for some generic existence results for constant $Q$-curvature metrics.
\end{remark}

Our second result demonstrates the existence of complete Ricci-flat metrics on $T(P)$.
\begin{theorem}\label{zero}
Let $V$ be a Fano K\"{a}hler-Einstein manifold with first Chern class $pa$, where $p\in\mathbb{Z}_+$, and $a\in H^2(V;\mathbb{Z})$ is an indivisible class. For $q=(q_1,q_2)\in\mathbb{Z}\oplus\mathbb{Z}$ with $|q_2|>0$, let $P_q$ be the principal $2$-torus bundle over $V$ with characteristic classes $(q_1a,q_2a)$. Then there exist a two-parameter family of complete Ricci-flat metrics on $T(P_q)$.
\end{theorem}

By the Bishop volume comparison theorem, a complete noncompact Ricci-flat manifold can have at most Euclidean volume growth. It can be shown that all the Ricci-flat manifolds in Theorem \ref{zero} have slower-than-Euclidean volume growth and quadratic curvature decay (cf. \cite{LotShe00}).

Now we turn to the associated $3$-sphere bundles. In general, given a principal $2$-torus bundle $P$, we can construct a $3$-sphere bundle $S(P)=P\times_{S^1\times S^1}S^3$ via the left action of $2$-torus $S^1\times S^1$ on unit $3$-sphere $S^3(1)=\{(z_1,z_2)\in\mathbb{C}^2:|z_1|^2+|z_2|^2=1\}$
\begin{align*}
(e^{i\alpha},e^{i\beta})\cdot(z_1,z_2)=(e^{i\alpha}z_1,e^{i\beta}z_2).
\end{align*}
Our third result demonstrates the existence of positive Einstein metrics on $S(P)$.

\begin{theorem}\label{pos}
Let $V$ be a Fano K\"{a}hler-Einstein manifold with first Chern class $pa$, where $p\in\mathbb{Z}_+$, and $a\in H^2(V;\mathbb{Z})$ is an indivisible class. For $q=(q_1,q_2)\in\mathbb{Z}\oplus\mathbb{Z}$ with $|q_1|>|q_2|>0$, let $P_q$ be the principal $2$-torus bundle over $V$ with characteristic classes $(q_1a,q_2a)$. Then there exists a smooth positive Einstein metric on $S(P_q)$.
\end{theorem}

\begin{remark}
In the case of $q_2=0$ and $0<|q_1|<p$, L\"{u}, Page and Pope \cite[\S3.3.2]{LvPagPop04} constructed a smooth positive Einstein metric on $S(P_q)$.
\end{remark}

\begin{remark}
A special case of Theorem \ref{pos}, where the base $V$ is $\mathbb{C}P^1$, was obtained by Hashimoto, Sakaguchi and Yasui \cite[Theorem 1]{HasSakYas05} in a different way.
\end{remark}

\begin{remark}
Let $\omega$ be a primitive $r$th root of unity, and let $s$ be an integer coprime to $r$. The cyclic group $\mathbb{Z}_r\subset S^1\times S^1$ of order $r$, generated by $(\omega,\omega^s)$, acts freely on $S^3$ from the right
\begin{align*}
(z_1,z_2)\cdot(\omega,\omega^s)=(z_1\omega,z_2\omega^s).
\end{align*}
It is clear from our construction that $\mathbb{Z}_r$ acts by isometry on the Einstein manifolds in Theorem \ref{pos}. The quotient manifolds, which are lens space bundles, inherit positive Einstein metrics.
\end{remark}

It is an interesting but difficult problem to study the moduli spaces of Einstein structures on the above $3$-sphere bundles. A first step in this direction is a diffeomorphism classification of the total spaces of these $3$-sphere bundles. By comparing the Einstein constants of Einstein metrics with unit volume on diffeomorphic total spaces, we may gain some information about the number of components of the moduli space .

The simplest case is that the base $V$ is $S^2$. There are only two $3$-sphere bundles over $S^2$ up to diffeomorphism \cite[\S26]{Ste99}. The total space of the trivial one is the product manifold $S^2\times S^3$ which is spin, while the total space of the nontrivial one is the twisted manifold $S^2\widetilde{\times}S^3$ which is non-spin. In Theorem \ref{pos}, when $q_1+q_2\equiv0$ mod $2$, the total space is spin, and hence is diffeomorphic to $S^2\times S^3$, while when $q_1+q_2\equiv1$ mod $2$, the total space is non-spin, and hence is diffeomorphic to $S^2\widetilde{\times}S^3$. It is well-known that the moduli space of Einstein structures on either manifold has infinitely many components \cite[\S2.2]{GibHarYas04}.

The situation becomes much more involved when the base $V$ has higher dimensions. We succeed in classifying the total spaces only when $V$ is $\mathbb{C}P^2$. In this case, the total space $W_q$ of the associated $3$-sphere bundle $S(P_q)$ is a simply-connected, closed $7$-manifold with $H^2(W_q;\mathbb{Z})\cong\mathbb{Z}$, $H^3(W_q;\mathbb{Z})=0$, and $H^4(W_q;\mathbb{Z})\cong\mathbb{Z}_{|q_1q_2|}$ generated by $u^2$, where $u$ is a generator of $H^2(W_q;\mathbb{Z})$. So $|q_1q_2|=|\widehat{q}_1\widehat{q}_2|$ is a necessary condition for $W_q$ to be homotopic to $W_{\widehat{q}}$. Furthermore, we can apply the classification result of Kreck and Stolz \cite{KreSto88, KreSto91} to show that

\begin{theorem}\label{cla}
Assume $K=q_1q_2$ and $L=q_1^2+q_2^2$.
\begin{enumerate}
  \item When $K=-\widehat{K}$, $W_q$ is homeomorphic (diffeomorphic) to $W_{\widehat{q}}$ iff $|K|=1$.
  \item When $K=\widehat{K}$, $W_q$ is homeomorphic (diffeomorphic) to $W_{\widehat{q}}$ iff $L\equiv\widehat{L}$ mod $2^{4-\mu(L)}\times3K$ (and $L+3^{\mu(L)}[\frac{L+1}{2}]^2\equiv\widehat{L}+3^{\mu(\widehat{L})}[\frac{\widehat{L}+1}{2}]^2$ mod $2^5\times3^{\mu(L)}\times7K$), where
\begin{align*}
\mu(L)=\left\{\begin{array}{ll}
0,&L\mbox{ is odd};\\
1,&L\mbox{ is even}.
\end{array}\right.
\end{align*}
\end{enumerate}
\end{theorem}

\begin{remark}
$W_q$ is not homotopic to any Aloff-Wallach space $M$. The argument goes as follows. Note that $M$ is spin, and $H^4(M;\mathbb{Z})$ is a finite cyclic group of odd order. For $W_q$ to be spin, $q_1+q_2$ has to be odd (see Lemma \ref{w}). So $q_1q_2$ is even, and $H^4(W_q;\mathbb{Z})$ is of even order. However, we do not know whether $W_q$ can be homeomorphic to any Eschenburg space. Note that the Kreck-Stolz invariants for a certain type of Eschenburg spaces can be found in \cite{Kru05}.
\end{remark}

Theorem \ref{cla} provides infinitely many pairs of homeomorphic manifolds which are not diffeomorphic, as well as infinitely many pairs of diffeomorphic manifolds. For instance,
\begin{example}[Spin case]\label{s}
For $r=48s+1$, $W_{(1,r(r+1))}$ and $W_{(r,r+1)}$ are homeomorphic. They are diffeomorphic iff $s\equiv0,3,4,6$ mod $7$.
\end{example}

\begin{example}[Non-spin case]\label{ns}
For $r=24s+1$, $W_{(2,2r(r+1))}$ and $W_{(2r,2(r+1))}$ are homeomorphic. They are diffeomorphic iff $s=4t$ with $t\equiv0,1,4,5,6$ mod $7$.
\end{example}

\begin{remark}\label{v}
Given a pair of diffeomorphic manifolds as in Examples \ref{s} and \ref{ns}, Theorem \ref{pos} asserts the existence of two positive Einstein metrics with unit volume on the underlying smooth manifold. But it is hard to verify whether they have the same Einstein constants or not (cf. Remark \ref{vv}). We believe that the Einstein constants should be different in general. Thus the moduli space of Einstein structures would have more than one component.
\end{remark}

The remainder of this paper is structured as follows. In \S\ref{pt}, we discuss a class of principal $2$-torus bundles over a Fano K\"{a}hler-Einstein manifold from both topological and geometric viewpoints. In \S\ref{fpt}, we compute Ricci curvatures of the warped product of an open interval and a principal $2$-torus bundle studied in \S\ref{pt}, and reduce the Einstein equation on the product to a system of ODEs. In \S\ref{e}, we find exact solutions to a subsystem of the Einstein system derived in \S\ref{fpt}. In \S\ref{nem}, we construct complete non-positive Einstein metrics on associated solid torus bundles. After that comes a proof of Theorem \ref{zq}. In \S\ref{pe}, we construct positive Einstein metrics on associated $3$-sphere bundles. We end the paper with a detailed calculation of characteristic classes of the $3$-sphere bundles and the $4$-ball bundles they bound, which leads to an argument for Theorem \ref{cla}.

\section{Principal $2$-torus bundles}\label{pt}

This section is devoted to a brief discussion of the topological and geometric properties of principal $2$-torus bundles involved in Theorems \ref{neg} and \ref{pos} (see \cite{WanZil90} for more details).

\subsection{Topology}

Let $(V^{2n},J,h)$ be a Fano K\"{a}hler-Einstein manifold of complex dimension $n>0$. Assume that the first Chern class $c_1(V,J)=pa$, where $p\in\mathbb{Z}_+$ is the Einstein constant of $h$, i.e., $\mbox{Ric}(h)=ph$, and $a\in H^2(V;\mathbb{Z})$ is an indivisible class. Notice that $H^2(V;\mathbb{Z})$ is torsion-free as $V$ is simply-connected.

Let $T^2$ be $2$-dimensional compact torus. We decompose once and for all $T^2$ as a product $S^1\times S^1$, and choose a basis $\{e_1,e_2\}$ for its Lie algebra $\mathfrak{t}^2$. Thus the set of isomorphism classes of principal $T^2$-bundles over $V$ is identified with $[V,\mathbb{C}P^\infty\times\mathbb{C}P^\infty]$, i.e., the set of homotopy classes of maps from $V$ to $\mathbb{C}P^\infty\times\mathbb{C}P^\infty$. Since $\mathbb{C}P^\infty\times\mathbb{C}P^\infty$ is the Eilenberg-MacLane space $K(\mathbb{Z}\oplus\mathbb{Z},2)$, $[V,\mathbb{C}P^\infty\times\mathbb{C}P^\infty]$ is isomorphic to $H^2(V;\mathbb{Z}\oplus\mathbb{Z})\cong H^2(V;\mathbb{Z})\oplus H^2(V;\mathbb{Z})$. Therefore a principal $T^2$-bundle over $V$ is classified by a pair of characteristic classes $(\upsilon_1,\upsilon_2)\in H^2(V;\mathbb{Z})\oplus H^2(V;\mathbb{Z})$.

The following describes how to construct the principal $T^2$-bundle classified by a given pair of characteristic classes. For $\upsilon_i\in H^2(V;\mathbb{Z})$, let $P_i$ be the principal $S^1$-bundle over $V$ with Euler class $\upsilon_i$. Their Cartesian product $P_1\times P_2$ is a principal $T^2$-bundle over $V\times V$. The pullback bundle $P$ of $P_1\times P_2$ by the diagonal map $V\rightarrow V\times V$ is the principal $T^2$-bundle over $V$ with characteristic classes $(\upsilon_1,\upsilon_2)$.

Observe that $T^2$ has a large automorphism group $GL(2,\mathbb{Z})$. This fact provides a simple way to construct new principal $T^2$-bundles out of old ones as follows. Given a principal $T^2$-bundle $P$ with characteristic classes $(\upsilon_1,\upsilon_2)$, we can change the $T^2$-action on $P$ via an element $A\in GL(2,\mathbb{Z})$, say $A=\left(\begin{array}{cc}
A_{11}&A_{12}\\
A_{21}&A_{22}
\end{array}\right)\in SL(2,\mathbb{Z})$.
\begin{align*}
\left.\begin{array}{cccc}
A:&S^1\times S^1&\rightarrow&S^1\times S^1\\
&(e^{i\alpha},e^{i\beta})&\mapsto&(e^{i(A_{11}\alpha+A_{12}\beta)},e^{i(A_{21}\alpha+A_{22}\beta)})
\end{array}\right.
\end{align*}
This yields a new principal $T^2$-bundle $\widetilde{P}$ with characteristic classes $(\widetilde{\upsilon}_1,\widetilde{\upsilon}_2)$, where $\widetilde{\upsilon}_j=(A^{-1})_{jk}\upsilon_k$. Notice that $P$ and $\widetilde{P}$ have the same manifold as their total spaces.

From now on, we will focus on a special class of principal $T^2$-bundles over $V$, whose characteristic classes are both integral multiples of $a$. Given $q=(q_1,q_2)\in\mathbb{Z}\oplus\mathbb{Z}\setminus\{(0,0)\}$, let $P_q$ be the principal $T^2$-bundle over $V$ with characteristic classes $(q_1a,q_2a)$. Denote by $q_0$ the greatest common divisor of $q_1$ and $q_2$, and let $P_{q_0}$ be the principal $S^1$-bundle over $V$ with Euler class $q_0a$. We have

\begin{proposition}\label{product}
The total spaces of $P_q$ and $P_{q_0}\times S^1$ are diffeomorphic to each other.
\end{proposition}

\begin{proof}
B\'{e}zout's identity asserts the existence of integers $r_1$ and $r_2$ such that $q_1r_1+q_2r_2=q_0$. Let
\begin{align*}
A=\left(\begin{array}{cc}
q_1/q_0&-r_2\\
q_2/q_0&r_1
\end{array}\right)
\in SL(2,\mathbb{Z}).
\end{align*}
Changing the $T^2$-action on $P_q$ via $A$ gives rise to a new principal $T^2$-bundle $\widetilde{P}$ over $V$ with characteristic classes $(q_0a,0)$. But $P_{q_0}\times S^1$ is also a principal $T^2$-bundle over $V$ with characteristic classes $(q_0a,0)$. Thus $\widetilde{P}$ is diffeomorphic to $P_{q_0}\times S^1$. The proposition follows from the fact that $P_q$ and $\widetilde{P}$ have the same total space.
\end{proof}

\subsection{Geometry}

Let $\widehat{\pi}_i:P_{q_i}\rightarrow V$ be the principal $S^1$-bundle with Euler class $q_ia$, $i=1,2$. There is a principal connection $\theta^i$ on $P_{q_i}$ with curvature form $\Omega_i=\widehat{\pi}_i^\ast\eta\otimes q_ie_i$, where $\eta$ is the K\"{a}hler form associated with $h$. If we denote by $\widetilde{e}_i$ the global vertical vector field on $P_{q_i}$ generated by $e_i$, i.e.,
\begin{eqnarray*}
\widetilde{e}_i(x)=\frac{d}{dt}|_{t=0}(x\cdot\exp(te_i)),&\forall x\in P_{q_i},
\end{eqnarray*}
then $\theta^i(\widetilde{e}_i)=e_i$.

Recall that the principal $T^2$-bundle $\widehat{\pi}:P_q\rightarrow V$ is the pullback bundle of the Cartesian product $\widehat{\pi}_1\times\widehat{\pi}_2:P_{q_1}\times P_{q_2}\rightarrow V\times V$ via the diagonal map $V\rightarrow V\times V$. Thus a principal connection $\theta$ on $P_q$ is given by the pullback of $\theta^1\times\theta^2$, and its curvature form is $\Omega=\widehat{\pi}^\ast\eta\otimes(q_1e_1+q_2e_2)$.

Let $B=(b_{ij})$ be an arbitrary positive-definite $2\times2$ symmetric matrix. It induces a left-invariant metric $\langle\cdot,\cdot\rangle_B$ on $T^2$ given by $\langle e_i,e_j\rangle_B=b_{ij}$. Consider now a bundle metric on $P_q$
\begin{eqnarray}\label{abm}
g=\langle\theta,\theta\rangle_B+c^2\widehat{\pi}^\ast h=\sum_{i,j=1}^2b_{ij}\theta^i\otimes\theta^j+c^2\widehat{\pi}^\ast h,
\end{eqnarray}
where $c$ is a positive constant, and the convention is $\theta^i(\widetilde{e}_j)=\delta^i_j$. Such a choice makes $\widehat{\pi}:(P_q,g)\rightarrow(V,c^2h)$ into a Riemannian submersion with totally geodesic fibers. More importantly, the curvature form $\Omega$ is parallel w.r.t. the base metric. Thus the principal connection $\theta$ satisfies the Yang-Mills condition. The relationship of the $(2,1)$ O'Neill tensor field $A$ (cf. \cite[9.20]{Bes87}) and the curvature form $\Omega$ becomes
\begin{align*}
A_XY=-\frac{1}{2}\Omega(X,Y)=-\frac{1}{2}\eta(\widehat{\pi}_\ast(X),\widehat{\pi}_\ast(Y))(q_1\widetilde{e}_1+q_2\widetilde{e}_2),
\end{align*}
where $X$ and $Y$ are arbitrary horizontal vector fields on $P_q$.

Now we compute the Ricci tensor of $(P_q,g)$. To do that we choose an orthonormal adapted basis $\{\check{E}_i\}_{1\le i\le2n}$ on $(V^{2n},J,h)$, i.e., $\check{E}_{n+j}=J\check{E}_j$, $1\le j\le n$, and let $E_i$ be the unique horizontal vector field on $P_q$ such that $\widehat{\pi}_\ast(E_i)=\check{E}_i$. This gives us a standard basis $\{\widetilde{e}_i,E_j\}_{1\le i\le2,1\le j\le2n}$ on $P_q$. In terms of such a basis, we have \cite[9.36]{Bes87}
\begin{lemma}\label{ric}
The non-vanishing components of the Ricci tensor of $(P_q,g)$ are
\begin{eqnarray*}
\mbox{Ric}(\widetilde{e}_i,\widetilde{e}_j)=\frac{n}{2c^4}U_iU_j,
&&\mbox{Ric}(E_k,E_l)=(p-\frac{\Delta}{2c^2})\delta_{kl},
\end{eqnarray*}
with $U=(U_1,U_2)=(q_1b_{11}+q_2b_{12}, q_1b_{12}+q_2b_{22})$ and $\Delta=q_1U_1+q_2U_2=q_1^2b_{11}+2q_1q_2b_{12}+q_2^2b_{22}$.
\end{lemma}
\begin{remark}\label{d}
$\Delta=qBq^T>0$ since $B$ is positive-definite, and $q=(q_1,q_2)\ne(0,0)$.
\end{remark}

Let $\widetilde{\mbox{Ric}}$ be the Ricci endomorphism of $TP_q$, i.e., $g(\widetilde{\mbox{Ric}}(X),Y)=\mbox{Ric}(X,Y)$, for $X,Y\in TP_q$. By Lemma \ref{ric}, we have
\begin{lemma}\label{ip}
\begin{align*}
\widetilde{\mbox{Ric}}(\widehat{e}_1)=\frac{n\Delta}{2c^4}\widehat{e}_1,&&\widetilde{\mbox{Ric}}(\widehat{e}_2)=0,&&\widetilde{\mbox{Ric}}(E_i)=(\frac{p}{c^2}-\frac{\Delta}{2c^4})E_i,
\end{align*}
where $\widehat{e}_1=\frac{1}{\sqrt{\Delta}}(q_1\widetilde{e}_1+q_2\widetilde{e}_2)$ and $\widehat{e}_2=\frac{1}{\sqrt{\alpha\Delta}}(U_2\widetilde{e}_1-U_1\widetilde{e}_2)$ are orthonormal vertical vector fields, and $\alpha=\det B=b_{11}b_{22}-b_{12}^2$. The scalar curvature of $g$ is
\begin{align*}
R(g)=\trace(\widetilde{\mbox{Ric}})=\frac{2n}{c^2}(p-\frac{\Delta}{4c^2}).
\end{align*}
\end{lemma}

Finally we derive a useful formula for $b_{ij}$ involving $\alpha$, $\Delta$ and $U_i$. Observe that the dual coframe of $\{\widehat{e}_i\}$ consists of $\widehat{\omega}^1=\frac{1}{\sqrt{\Delta}}(U_1\theta^1+U_2\theta^2)$ and $\widehat{\omega}^2=\sqrt{\frac{\alpha}{\Delta}}(q_2\theta^1-q_1\theta^2)$, i.e., $\widehat{\omega}^i(\widehat{e}_j)=\delta^i_j$. In terms of $\{\widehat{\omega}^i\}$, the vertical components of the metric $g$ \eqref{abm} takes a diagonal form
\begin{align*}
\sum_{i,j=1}^2b_{ij}\theta^i\otimes\theta^j=\widehat{\omega}^1\otimes\widehat{\omega}^1+\widehat{\omega}^2\otimes\widehat{\omega}^2.
\end{align*}
Comparing the coefficients on both sides yields
\begin{lemma}\label{b}
\begin{align*}
b_{11}=\frac{U_1^2+q_2^2\alpha}{\Delta},&&b_{12}=\frac{U_1U_2-q_1q_2\alpha}{\Delta},
&&b_{22}=\frac{U_2^2+q_1^2\alpha}{\Delta}.
\end{align*}
\end{lemma}

\section{A family of principal $T^2$-bundles}\label{fpt}

Let $P_q$ be a principal $T^2$-bundle as in \S\ref{pt}, and let $T(P_q)$ and $S(P_q)$ be respectively the associated solid torus bundle and $3$-sphere bundle as in \S\ref{intr}. Notice that there exist diffeomorphisms
\begin{align*}
&T_0(P_q)=P_q\times_{S^1\times S^1}(B^2\times S^1\backslash\{\{0\}\times S^1\})\cong I\times P_q,\\
&S_0(P_q)=P_q\times_{S^1\times S^1}(S^3\backslash\{\{0\}\times S^1,S^1\times\{0\}\})\cong I\times P_q,
\end{align*}
where $I\subset\mathbb{R}$ is an open interval. Thus our strategy for proving Theorems \ref{neg}, \ref{zero} and \ref{pos} consists of constructing local metrics on $T_0(P_q)$ and $S_0(P_q)$, and then extending them to be smooth Einstein metrics on $T(P_q)$ and $S(P_q)$ respectively.

In this section, we study the geometry of product $\widehat{M}=I\times P_q$ endowed with metric $\widehat{g}=dt^2+g_t$, where
\begin{eqnarray}\label{bm}
g_t=\sum_{i,j=1}^2b_{ij}(t)\theta^i\otimes\theta^j+c(t)^2\widehat{\pi}^\ast h,&&t\in I,
\end{eqnarray}
$B(t)=(b_{ij}(t))$ is a smooth one-parameter family of positive-definite $2\times2$ symmetric matrices, and $c(t)$ is a smooth positive function of $t$. We will reduce the Einstein equation \eqref{ee} on $(\widehat{M},\widehat{g})$ to a system of ODEs.

Denote by $L_t$ the self-adjoint shape operator of hypersurfaces $\Sigma_t=\{t\}\times P_q$. By definition,
\begin{eqnarray*}
L_tX=\widehat{\nabla}_XN,&X\in T\Sigma_t,
\end{eqnarray*}
where $\widehat{\nabla}$ is the Levi-Civita connection of $\widehat{g}$, and $N=\frac{\partial}{\partial t}$ is the unit normal vector.
\begin{lemma}\label{so}
For $X,Y\in T\Sigma_t$, we have
\begin{eqnarray*}
g_t(L_tX,Y)=\frac{1}{2}g_t'(X,Y),
\end{eqnarray*}
where $'$ denotes $\frac{\partial}{\partial t}$.
\end{lemma}
\begin{proof}
Given a vector field $Z$ on $\Sigma_t$, we denote by $\widehat{Z}$ its unique $t$-independent extension to $\widehat{M}$, such that the Lie bracket $[\widehat{Z},N]=0$, i.e., $\widehat{\nabla}_N\widehat{Z}=\widehat{\nabla}_{\widehat{Z}}N$. Thus
\begin{eqnarray*}
g_t'(X,Y)&=&\frac{\partial}{\partial t}\widehat{g}(\widehat{X},\widehat{Y})\\
&=&\widehat{g}(\widehat{\nabla}_N\widehat{X},\widehat{Y})+\widehat{g}(\widehat{X},\widehat{\nabla}_N\widehat{Y})\\
&=&\widehat{g}(\widehat{\nabla}_{\widehat{X}}N,\widehat{Y})+\widehat{g}(\widehat{X},\widehat{\nabla}_{\widehat{Y}}N)\\
&=&g_t(L_tX,Y)+g_t(X,L_tY)\\
&=&2g_t(L_tX,Y).
\end{eqnarray*}
\end{proof}
\begin{lemma}\label{sh}
On the standard basis $\{\widetilde{e}_i,E_j\}_{1\le i\le2,1\le j\le2n}$, $L_t$ can be written in a matrix form
\begin{align*}
(L_t\widetilde{e}_1,L_t\widetilde{e}_2,L_tE_1,\cdots,L_tE_{2n})=(\widetilde{e}_1,\widetilde{e}_2,E_1,\cdots,E_{2n})\left(\begin{array}{cc}
\frac{1}{2}\Psi&0\\
0&\frac{c'}{c}\mbox{I}_{2n}\end{array}
\right)^T,
\end{align*}
where $\Psi=B'B^{-1}$, $B^{-1}=(b^{ij})$ is the inverse matrix of $B$, and $\mbox{I}_{2n}$ is the $2n\times2n$ identity matrix. In particular,
\begin{eqnarray*}
\trace(L_t)=\frac{1}{2}\trace(\Psi)+2n\frac{c'}{c},&&
\trace(L_t^2)=\frac{1}{4}\trace(\Psi^2)+2n\frac{c'^2}{c^2}.
\end{eqnarray*}

\end{lemma}
\begin{proof}
It follows from Lemma \ref{so} and \eqref{bm} that on $T\Sigma_t$
\begin{align}\label{sod}
g_t(L_t(\cdot),\cdot)=\frac{1}{2}\sum_{i,j=1}^2b_{ij}'\theta^i\otimes\theta^j+cc'\widehat{\pi}^\ast h.
\end{align}
Clearly, $g_t(L_t\widetilde{e}_i,E_j)=0$. We need to determine the remaining components.

First, assume that $L_t\widetilde{e}_i=\sum_{j=1}^2\lambda_{ij}\widetilde{e}_j$, $i=1,2$. By \eqref{bm} and \eqref{sod}, we have $\sum_{k=1}^2\lambda_{ik}b_{kj}=\frac{1}{2}b_{ij}'$. Thus $\lambda_{ij}=\frac{1}{2}\sum_{k=1}^2b_{ik}'b^{kj}$.

Second, assume that $L_tE_i=\sum_{j=1}^{2n}\mu_{ij}E_j$, $1\le i\le2n$. By \eqref{bm} and \eqref{sod}, we have $c^2\mu_{ij}=cc'\delta_{ij}$. Thus $\mu_{ij}=\frac{c'}{c}\delta_{ij}$.
\end{proof}

\begin{corollary}\label{trlp}
On the standard basis $\{\widetilde{e}_i,E_j\}_{1\le i\le2,1\le j\le 2n}$, $L_t'$ can be written in a matrix form
\begin{align*}
(L_t'\widetilde{e}_1,L_t'\widetilde{e}_2,L_t'E_1,\cdots,L_t'E_{2n})=(\widetilde{e}_1,\widetilde{e}_2,E_1,\cdots,E_{2n})\left(\begin{array}{cc}
\frac{1}{2}\Psi'&0\\
0&(\frac{c''}{c}-\frac{c'^2}{c^2})\mbox{I}_{2n}
\end{array}\right)^T.
\end{align*}
In particular,
\begin{align*}
\trace(L_t')=\frac{1}{2}\trace(\Psi')+2n(\frac{c''}{c}-\frac{c'^2}{c^2}).
\end{align*}
\end{corollary}

Now we compute the Ricci curvatures of $(\widehat{M},\widehat{g})$. For vector fields $X,Y,Z\in T\Sigma_t$, we have
\begin{eqnarray*}
\widehat{g}(\widehat{R}(X,Y)N,Z)&=&\widehat{g}(\widehat{\nabla}_Y\widehat{\nabla}_XN-\widehat{\nabla}_X\widehat{\nabla}_YN
+\widehat{\nabla}_{[X,Y]}N,Z)\\
&=&\widehat{g}(\widehat{\nabla}_YL_tX-\widehat{\nabla}_XL_tY+L_t[X,Y],Z)\\
&=&g_t(\nabla^t_YL_tX-\nabla^t_XL_tY+L_t[X,Y],Z),
\end{eqnarray*}
where $\nabla^t$ is the Levi-Civita connection of $g_t$. Thus
\begin{eqnarray*}
\widehat{\mbox{Ric}}(Z,N)&=&\sum_{i,j=1}^2b^{ij}\widehat{g}(\widehat{R}(Z,\widetilde{e}_i)N,\widetilde{e}_j)+\sum_{k=1}^{2n}c^{-2}\widehat{g}(\widehat{R}(Z,E_k)N,E_k)\\
&=&\sum_{i,j=1}^2b^{ij}g_t(\nabla^t_{\widetilde{e}_i}L_tZ-\nabla^t_ZL_t\widetilde{e}_i+L_t[Z,\widetilde{e}_i],\widetilde{e}_j)\\
&&+\sum_{k=1}^{2n}c^{-2}g_t(\nabla^t_{E_k}L_tZ-\nabla^t_ZL_tE_k+L_t[Z,E_k],E_k).
\end{eqnarray*}

\begin{lemma}
$\widehat{\mbox{Ric}}(\widetilde{e}_l,N)=0$, $l=1,2$.
\end{lemma}

\begin{proof}
Notice that
\begin{eqnarray*}
\widehat{\mbox{Ric}}(\widetilde{e}_l,N)=\sum_{i,j=1}^2b^{ij}g_t(\nabla^t_{\widetilde{e}_i}L_t\widetilde{e}_l-\nabla^t_{\widetilde{e}_l}L_t\widetilde{e}_i,\widetilde{e}_j)+\sum_{k=1}^{2n}c^{-2}g_t(\nabla^t_{E_k}L_t\widetilde{e}_l-\nabla^t_{\widetilde{e}_l}L_tE_k+L_t[\widetilde{e}_l,E_k],E_k),
\end{eqnarray*}
where we have used the fact $[\widetilde{e}_l,\widetilde{e}_i]=0$ since the fiber of the Riemannian submersion $\widehat{\pi}:(\Sigma_t,g_t)\rightarrow(V,c^2h)$ is an abelian group. But $[\widetilde{e}_l,E_k]$ is vertical, so is $L_t[\widetilde{e}_l,E_k]$. Thus $g_t(L_t[\widetilde{e}_l,E_k],E_k)=0$.
\begin{claim}\label{ab}
$g_t(\nabla^t_{\widetilde{e}_i}\widetilde{e}_r,\widetilde{e}_j)=0$.
\end{claim}
\begin{proof}
In fact,
\begin{eqnarray*}
g_t(\nabla^t_{\widetilde{e}_i}\widetilde{e}_r,\widetilde{e}_j)&=&-g_t(\widetilde{e}_r,\nabla^t_{\widetilde{e}_i}\widetilde{e}_j)\\
&=&-g_t(\widetilde{e}_r,\nabla^t_{\widetilde{e}_j}\widetilde{e}_i)\\
&=&g_t(\nabla^t_{\widetilde{e}_j}\widetilde{e}_r,\widetilde{e}_i)\\
&=&g_t(\nabla^t_{\widetilde{e}_r}\widetilde{e}_j,\widetilde{e}_i)\\
&=&-g_t(\widetilde{e}_j,\nabla^t_{\widetilde{e}_r}\widetilde{e}_i)\\
&=&-g_t(\widetilde{e}_j,\nabla^t_{\widetilde{e}_i}\widetilde{e}_r).
\end{eqnarray*}
\end{proof}
It follows from Claim \ref{ab} that
\begin{eqnarray*}
g_t(\nabla^t_{\widetilde{e}_i}L_t\widetilde{e}_l,\widetilde{e}_j)=\sum_{r=1}^2\lambda_{lr}g_t(\nabla^t_{\widetilde{e}_i}\widetilde{e}_r,\widetilde{e}_j)=0.
\end{eqnarray*}
Similarly, we have $g_t(\nabla^t_{\widetilde{e}_l}L_t\widetilde{e}_i,\widetilde{e}_j)=0$. Moreover,
\begin{align*}
&g_t(\nabla^t_{E_k}L_t\widetilde{e}_l,E_k)=g_t(A^t_{E_k}L_t\widetilde{e}_l,E_k)=-g_t(A^t_{E_k}E_k,L_t\widetilde{e}_l)=0,\\
&g_t(\nabla^t_{\widetilde{e}_l}L_tE_k,E_k)=\frac{c'}{c}g_t(\nabla^t_{\widetilde{e}_l}E_k,E_k)=\frac{c'}{2c}\widetilde{e}_l(g_t(E_k,E_k))=0,
\end{align*}
where $A^t$ is the $(2,1)$ O'Neill tensor field associated with $g_t$ (cf. \cite[9.20]{Bes87}).

To sum up, we see that $\widehat{\mbox{Ric}}(\widetilde{e}_l,N)=0$.
\end{proof}

\begin{lemma}
$\widehat{\mbox{Ric}}(E_l,N)=0$, $1\le l\le2n$.
\end{lemma}

\begin{proof}
We have
\begin{eqnarray*}
\widehat{\mbox{Ric}}(E_l,N)&=&\sum_{i,j=1}^2b^{ij}g_t(\nabla^t_{\widetilde{e}_i}L_tE_l-\nabla^t_{E_l}L_t\widetilde{e}_i+L_t[E_l,\widetilde{e}_i],\widetilde{e}_j)\\
&&+\sum_{k=1}^{2n}c^{-2}g_t(\nabla^t_{E_k}L_tE_l-\nabla^t_{E_l}L_tE_k+L_t[E_l,E_k],E_k).
\end{eqnarray*}
Observe that $L_tE_l$ is horizontal, so is $\nabla^t_{\widetilde{e}_i}L_tE_l$ as the fiber is totally geodesic (cf. \cite[(9.25b)]{Bes87}). Thus $g_t(\nabla^t_{\widetilde{e}_i}L_tE_l,\widetilde{e}_j)=0$. Moreover,
\begin{eqnarray*}
\sum_{i,j=1}^2b^{ij}g_t(\nabla^t_{E_l}L_t\widetilde{e}_i,\widetilde{e}_j)&=&\sum_{i,j,r=1}^2b^{ij}\lambda_{ir}g_t(\nabla^t_{E_l}\widetilde{e}_r,\widetilde{e}_j)\\
&=&-\sum_{i,j,r=1}^2b^{ij}\lambda_{ir}g_t(\widetilde{e}_r,\nabla^t_{E_l}\widetilde{e}_j)\\
&=&-\sum_{i,j,r=1}^2b^{ir}\lambda_{ij}g_t(\widetilde{e}_j,\nabla^t_{E_l}\widetilde{e}_r).
\end{eqnarray*}
Hence
\begin{eqnarray*}
0&=&\sum_{i,j,r=1}^2(b^{ij}\lambda_{ir}+b^{ir}\lambda_{ij})g_t(\widetilde{e}_j,\nabla^t_{E_l}\widetilde{e}_r)\\
&=&\sum_{i,j,k,r=1}^2\frac{1}{2}(b_{ik}'b^{kr}b^{ij}+b_{ik}'b^{kj}b^{ir})g_t(\widetilde{e}_j,\nabla^t_{E_l}\widetilde{e}_r)\\
&=&\sum_{i,j,k,r=1}^2b_{ik}'b^{kr}b^{ij}g_t(\widetilde{e}_j,\nabla^t_{E_l}\widetilde{e}_r)\\
&=&\sum_{i,j,r=1}^22\lambda_{ir}b^{ij}g_t(\widetilde{e}_j,\nabla^t_{E_l}\widetilde{e}_r).
\end{eqnarray*}
It follows that
\begin{align*}
\sum_{i,j=1}^2b^{ij}g_t(\nabla^t_{E_l}L_t\widetilde{e}_i,\widetilde{e}_j)=0.
\end{align*}
On the other hand, since $\nabla^t_{\widetilde{e}_i}E_l$ is horizontal, and $L_t\widetilde{e}_j$ is vertical, we have
\begin{align*}
g_t(L_t[E_l,\widetilde{e}_i],\widetilde{e}_j)=g_t([E_l,\widetilde{e}_i],L_t\widetilde{e}_j)=g_t(\nabla^t_{E_l}\widetilde{e}_i,L_t\widetilde{e}_j)=-g_t(\widetilde{e}_i,\nabla^t_{E_l}L_t\widetilde{e}_j).
\end{align*}
Therefore,
\begin{align*}
\sum_{i,j=1}^2b^{ij}g_t(L_t[E_l,\widetilde{e}_i],\widetilde{e}_j)=-\sum_{i,j=1}^2b^{ij}g_t(\widetilde{e}_i,\nabla^t_{E_l}L_t\widetilde{e}_j)=0.
\end{align*}
Finally, we have
\begin{align*}
g_t(\nabla^t_{E_k}L_tE_l-\nabla^t_{E_l}L_tE_k+L_t[E_l,E_k],E_k)
=\frac{c'}{c}g_t(\nabla^t_{E_k}E_l-\nabla^t_{E_l}E_k+[E_l,E_k],E_k)=0.
\end{align*}

To sum up, we see that $\widehat{\mbox{Ric}}(E_l,N)=0$.
\end{proof}

By Lemmas \ref{ric}, \ref{sh} and Corollary \ref{trlp}, it is straightforward to show that \cite[\S2]{EscWan00}

\begin{lemma}
The non-vanishing components of the Ricci tensor of $(\widehat{M},\widehat{g})$ are
\begin{align}
&\widehat{\mbox{Ric}}(N,N)=-\frac{1}{2}\trace(\Psi')-\frac{1}{4}\trace(\Psi^2)-2n\frac{c''}{c},\\
&\widehat{\mbox{Ric}}(\widetilde{e}_i,\widetilde{e}_j)=\frac{n}{2c^4}U_iU_j-\frac{1}{2}b_{ij}'(\frac{1}{2}\trace(\Psi)
        +2n\frac{c'}{c})-\frac{1}{2}b_{ij}''+\frac{1}{2}\sum_{k,l=1}^2b_{ik}'b_{jl}'b^{kl},\\
&\widehat{\mbox{Ric}}(E_k,E_l)=(p-\frac{\Delta}{2c^2}-cc'(\frac{1}{2}\trace(\Psi)
        +2n\frac{c'}{c})-cc''+c'^2)\delta_{kl}.
\end{align}
\end{lemma}

Thus the Einstein equation $\widehat{\mbox{Ric}}=\epsilon\widehat{g}$ on $(\widehat{M}^{2n+3},\widehat{g})$ reduces to a system of ODEs
\begin{eqnarray}
\label{phh}&&-\frac{1}{2}\trace(\Psi')-\frac{1}{4}\trace(\Psi^2)-2n\frac{c''}{c}=\epsilon,\\
\label{pse}&&-\frac{1}{2}\Psi'-\frac{1}{2}(\frac{1}{2}\trace(\Psi)+2n\frac{c'}{c})\Psi+\frac{n}{2c^4}U^Tq=\epsilon\mbox{I}_2,\\
\label{pba}&&-(\frac{1}{2}\trace(\Psi)+2n\frac{c'}{c})\frac{c'}{c}-\frac{c''}{c}+\frac{c'^2}{c^2}+\frac{p}{c^2}
-\frac{\Delta}{2c^4}=\epsilon.
\end{eqnarray}

We close this section by quoting a couple of elementary identities involving $\Psi$ and $\alpha$.
\begin{proposition}\label{ei}
\begin{eqnarray*}
\trace(\Psi)=(\log\alpha)',&&\trace(\Psi^2)=\trace(\Psi)^2-2\det\Psi.
\end{eqnarray*}
\end{proposition}

\section{Exact solutions of the Einstein system}\label{e}

In this section, we construct exact solutions of the Einstein system \eqref{phh}-\eqref{pba}. To do that we introduce a new coordinate $s$ defined by $ds=\sqrt{\alpha}dt$, and let $\beta(s)=c(t)^2$. Then the metric $\widehat{g}$ takes the form
\begin{align}\label{em}
\widehat{g}=\alpha(s)^{-1}ds^2+\sum_{i,j=1}^2b_{ij}(s)\theta^i\otimes\theta^j+\beta(s)\widehat{\pi}^\ast h,
\end{align}
and the Einstein system \eqref{phh}-\eqref{pba} becomes
\begin{align}
\label{hh}&&n\alpha(-\frac{\ddot{\beta}}{\beta}+\frac{1}{2}(\frac{\dot{\beta}}{\beta})^2)
-n\frac{\dot{\alpha}}{2}\frac{\dot{\beta}}{\beta}-\frac{\ddot{\alpha}}{2}+\frac{\Upsilon}{2}=\epsilon,\\
\label{se}&&\frac{\alpha}{2}(-n\frac{\dot{\beta}}{\beta}\Phi-\dot{\Phi})-\frac{\dot{\alpha}}{2}\Phi
+\frac{n}{2\beta^2}U^Tq=\epsilon\mbox{I}_2,\\
\label{ba}&&\frac{\alpha}{2}(-\frac{\ddot{\beta}}{\beta}-(n-1)(\frac{\dot{\beta}}{\beta})^2)
-\frac{\dot{\alpha}}{2}\frac{\dot{\beta}}{\beta}+\frac{p}{\beta}-\frac{\Delta}{2\beta^2}=\epsilon,
\end{align}
where $^.$ denotes $\frac{d}{ds}$, $\Upsilon=\det\dot{B}$ and $\Phi=\dot{B}B^{-1}$. Here we have used Proposition \ref{ei}.

Thanks to Back (see \cite[Lemma 2.4]{EscWan00}) we can restrict ourselves to the subsystem consisting of \eqref{se}-\eqref{ba}, and look for only those solutions which can be extended to at least one endpoint of $I$. However, to the best of our knowledge, there is no general way to explicitly solve such a coupled non-linear system. Fortunately our work in \cite{Che09} (cf. Remark \ref{pw}) sheds some light on the problem. Motivated by \cite[Example 4.6.4]{Che09} our starting point is to assume $\beta$ is a linear function of $s$.
\begin{assumption}\label{as}
$\beta=\kappa s$ for some positive constant $\kappa$.
\end{assumption}
Under this assumption, we are able to find exact solutions of \eqref{se}-\eqref{ba} in the following way. First we take the trace of both sides of \eqref{se} to get
\begin{align}\label{tr}
-n\frac{\dot{\alpha}}{2}\frac{\dot{\beta}}{\beta}-\frac{\ddot{\alpha}}{2}+n\frac{\Delta}{2\beta^2}=2\epsilon.
\end{align}
Then $n\times\eqref{ba}+\eqref{tr}$ yields
\begin{align}\label{2tc}
\ddot{\alpha}=-2n\dot{\alpha}\frac{\dot{\beta}}{\beta}
+n\alpha(-\frac{\ddot{\beta}}{\beta}-(n-1)(\frac{\dot{\beta}}{\beta})^2)+\frac{2np}{\beta}-2\epsilon(n+2).
\end{align}
By Assumption \ref{as}, \eqref{2tc} becomes
\begin{align}\label{al}
\ddot{\alpha}=-\frac{2n}{s}\dot{\alpha}-\frac{n(n-1)}{s^2}\alpha+\frac{2np}{\kappa}s^{-1}-2\epsilon(n+2).
\end{align}
It is not hard to see that \eqref{al} has the general solution
\begin{eqnarray}\label{als}
\alpha=-\frac{2\epsilon}{n+1}s^2+\frac{2p}{\kappa(n+1)}s+c_1s^{1-n}+c_2s^{-n},
\end{eqnarray}
where $c_i$'s are integration constants. Substituting \eqref{als} into \eqref{ba} gives
\begin{align}\label{de}
\Delta=\frac{2p\kappa}{n+1}s+\kappa^2c_2s^{-n}.
\end{align}
Now letting $\widehat{q}=(q_2,-q_1)$, and multiplying $\widehat{q}^T$ on both sides of \eqref{se}, we get
\begin{align*}
(n\alpha\frac{\dot{\beta}}{\beta}\Phi+\alpha\dot{\Phi}+\dot{\alpha}\Phi)\widehat{q}^T=-2\epsilon\widehat{q}^T.
\end{align*}
i.e.,
\begin{align*}
\frac{d}{ds}(\alpha\beta^n\Phi\widehat{q}^T)=-2\epsilon\beta^n\widehat{q}^T.
\end{align*}
Thus
\begin{align*}
\alpha\beta^n\Phi\widehat{q}^T=-2\epsilon\int\beta^nds\cdot\widehat{q}^T+w^T
\end{align*}
for some constant vector $w=(w_1,w_2)$. In details, we have
\begin{eqnarray}
\label{2tv1}&&\dot{b_{11}}U_2-\dot{b_{12}}U_1=v_1,\\
\label{2tv2}&&\dot{b_{12}}U_2-\dot{b_{22}}U_1=v_2,
\end{eqnarray}
where
\begin{align*}
&v_1=\beta^{-n}(-2q_2\epsilon\int\beta^nds+w_1)=-\frac{2q_2\epsilon}{n+1}s+w_1\kappa^{-n}s^{-n},\\
&v_2=\beta^{-n}(2q_1\epsilon\int\beta^nds+w_2)=\frac{2q_1\epsilon}{n+1}s+w_2\kappa^{-n}s^{-n}.
\end{align*}

\begin{remark}
In the light of Lemma \ref{b}, it remains to determine $U_1$ and $U_2$ from \eqref{2tv1} and \eqref{2tv2}, and to verify that they satisfy the identity $q_1U_1+q_2U_2=\Delta$ (cf. Remark \ref{d}). After this is done, we can define $B=(b_{ij})$ as in Lemma \ref{b}, whose determinant turns out to be $\alpha$ automatically.
\end{remark}

Notice that $q_1\cdot\eqref{2tv1}+q_2\cdot\eqref{2tv2}$ yields
\begin{align}
\label{2tcs}\dot{U_1}U_2-U_1\dot{U_2}=\psi\beta^{-n},
\end{align}
where $\psi=q_1w_1+q_2w_2$ is a constant. We divide our discussion into two cases $\psi\ne0$ and $\psi=0$. The latter case is postponed until we place smooth collapse conditions on $b_{ij}$ in \S\ref{nem}.

Now assume $\psi\ne0$. It follows from Lemma \ref{b} that
\begin{eqnarray*}
&&\dot{b_{11}}=\frac{2U_1\dot{U_1}+q_2^2\dot{\alpha}}{\Delta}-\frac{U_1^2+q_2^2\alpha}{\Delta^2}\dot{\Delta},\\
&&\dot{b_{22}}=\frac{2U_2\dot{U_2}+q_1^2\dot{\alpha}}{\Delta}-\frac{U_2^2+q_1^2\alpha}{\Delta^2}\dot{\Delta},\\
&&\dot{b_{12}}=\frac{\dot{U_1}U_2+U_1\dot{U_2}-q_1q_2\dot{\alpha}}{\Delta}
-\frac{U_1U_2-q_1q_2\alpha}{\Delta^2}\dot{\Delta}.
\end{eqnarray*}
Substituting them into \eqref{2tv1} and \eqref{2tv2} yields
\begin{eqnarray}
\label{2tu}&&U_1=\psi^{-1}\beta^n(v_1\Delta-q_2(\dot{\alpha}\Delta-\alpha\dot{\Delta})),\\
\label{2tv}&&U_2=\psi^{-1}\beta^n(v_2\Delta+q_1(\dot{\alpha}\Delta-\alpha\dot{\Delta})),
\end{eqnarray}
provided \eqref{2tcs} is valid. Plugging them back into \eqref{2tcs} leads to the following consistency condition
\begin{eqnarray}
\label{2tccs}\psi^2=\beta^{2n}\dot{\Delta}(\dot{\alpha}\Delta-\alpha\dot{\Delta})
-\beta^n\Delta(\beta^n(\dot{\alpha}\Delta-\alpha\dot{\Delta}))\dot{}-2\epsilon\beta^{2n}\Delta^2.
\end{eqnarray}
The right-hand side of \eqref{2tccs} looks quite complicated, while the left-hand side is just a constant. It is thus surprising to see that when we plug \eqref{als} and \eqref{de} into \eqref{2tccs}, the consistency condition reads
\begin{eqnarray}
\label{2tcss}\psi^2=2(n+1)\kappa^{2n+3}c_2(pc_1+c_2\epsilon\kappa)>0
\end{eqnarray}
as $\psi\ne0$. Assuming this holds true, we can explicitly write down $U_1$ and $U_2$ respectively according to \eqref{2tu} and \eqref{2tv}. The upshot is
\begin{align}
\label{2tue}&U_1=\frac{2\kappa}{\psi(n+1)}(pw_1+npq_2\kappa^nc_1+q_2\epsilon\kappa^{n+1}c_2(n+1))s
+\frac{\kappa^2c_2}{\psi}(w_1-q_2\kappa^nc_1)s^{-n},\\
\label{2tve}&U_2=\frac{2\kappa}{\psi(n+1)}(pw_2-npq_1\kappa^nc_1-q_1\epsilon\kappa^{n+1}c_2(n+1))s
+\frac{\kappa^2c_2}{\psi}(w_2+q_1\kappa^nc_1)s^{-n}.
\end{align}
It is clear that $q_1U_1+q_2U_2=\Delta$ (cf. \eqref{de}).

\begin{proposition}\label{ls}
\eqref{als}, \eqref{de}, \eqref{2tue}, and \eqref{2tve}, together with the consistency condition \eqref{2tcss}, give rise to an exact solution of \eqref{se} and \eqref{ba} if we define $b_{ij}$ as in Lemma \ref{b}.
\end{proposition}

\section{Non-positive Einstein manifolds}\label{nem}

The data given in Proposition \ref{ls} yield local metrics on $I\times P_q$. In general, they are not Einstein metrics since the solutions may fail to satisfy \eqref{hh}. As noted in the paragraph preceding Assumption \ref{as}, one way to pick out Einstein metrics from these candidates is to smoothly extend them to one endpoint of $I$ by choosing parameters appropriately. This amounts to collapse a circle in each fiber of $P_q$. In practice, we collapse the circle in $2$-torus corresponding to characteristic class $q_1a$ at the left endpoint of $I$, say $s_1\ge0$. This gives rise to the following smoothness conditions on $b_{ij}(t)$:
\begin{enumerate}
  \item $b_{ij}(t)$'s are smooth and even in $t$ around $t_1=t(s_1)$;
  \item $b_{11}(t_1)=b_{12}(t_1)=0$, $b_{22}(t_1)>0$, and $\frac{d}{dt}|_{t=t_1}\sqrt{b_{11}(t)}=1$.
\end{enumerate}
The smoothness condition on $c(t)$ will be given in a while.
\begin{remark}\label{hhh}
Condition $(2)$ implies $\alpha(s_1)=U_1(s_1)=0$.
\end{remark}

The rest of this section is devoted to constructing complete non-positive Einstein manifolds. In particular, we assume $\epsilon\le0$. Our discussion is again divided into two cases $\psi\ne0$ and $\psi=0$.

\subsection{$\psi\ne0$}\label{pn}

It follows from \eqref{als} that
\begin{align}\label{aq}
\alpha=Qs^{-n},
\end{align}
where
\begin{align*}
Q=-\frac{2\epsilon}{n+1}s^{n+2}+\frac{2p}{\kappa(n+1)}s^{n+1}+c_1s+c_2.
\end{align*}
By \eqref{2tcss}, $Q(0)=c_2\ne0$. Thus the identity $\alpha(s_1)=0$ forces $s_1>0$ and $Q(s_1)=0$. It follows that
\begin{align}\label{c2}
c_2=\frac{2\epsilon}{n+1}s_1^{n+2}-\frac{2p}{\kappa(n+1)}s_1^{n+1}-c_1s_1.
\end{align}
\begin{remark}
Later on, we will see that $s_1$ is actually the largest positive root of $Q$.
\end{remark}

Notice that $c(t_1)>0$ as $\beta(s_1)=\kappa s_1>0$. Thus the smoothness condition on $c(t)$ is
\begin{itemize}
  \item $c(t)$ is smooth and even in $t$ around $t_1$.
\end{itemize}

On the other hand, it follows from \eqref{2tue} that $U_1(s_1)=0$ iff
\begin{eqnarray}\label{u1star}
\frac{2}{n+1}(pw_1+npq_2\kappa^nc_1+\epsilon q_2\kappa^{n+1}c_2(n+1))s_1^{n+1}
=-\kappa c_2(w_1-q_2\kappa^nc_1).
\end{eqnarray}
Using \eqref{c2} to eliminate $c_2$ in \eqref{u1star} yields
\begin{align}\label{ec2}
 -\kappa s_1(c_1-\frac{2\epsilon}{n+1}s_1^{n+1})(w_1-q_2\kappa^{n-1}(\kappa c_1+2ps_1^n-2\epsilon\kappa s_1^{n+1}))=0.
\end{align}
\begin{claim}\label{cc1}
$c_1-\frac{2\epsilon}{n+1}s_1^{n+1}\ne0$.
\end{claim}
\begin{proof}
If $c_1=\frac{2\epsilon}{n+1}s_1^{n+1}$, then it follows from \eqref{c2} that $c_2=-\frac{2p}{\kappa(n+1)}s_1^{n+1}$. So $pc_1+\epsilon\kappa c_2=0$. This contradicts \eqref{2tcss}.
\end{proof}
By Claim \ref{cc1}, we get from \eqref{ec2} that
\begin{align}\label{w1}
w_1=q_2\kappa^{n-1}(\kappa c_1+2ps_1^n-2\epsilon\kappa s_1^{n+1}).
\end{align}
Now it follows from \eqref{2tve}, \eqref{c2}, \eqref{w1}, and the fact $w_2=q_2^{-1}(\psi-q_1w_1)$ that
\begin{align*}
U_2(s_1)=-q_2^{-1}s_1^{1-n}\kappa^2(c_1-\frac{2\epsilon}{n+1}s_1^{n+1}).
\end{align*}
Thus
\begin{align*}
\Delta(s_1)=q_2U_2(s_1)=-s_1^{1-n}\kappa^2(c_1-\frac{2\epsilon}{n+1}s_1^{n+1})>0.
\end{align*}
So $c_1-\frac{2\epsilon}{n+1}s_1^{n+1}<0$, i.e.,
\begin{align}\label{c1n}
c_1<\frac{2\epsilon}{n+1}s_1^{n+1}\le0.
\end{align}
\begin{remark}
This upper bound on $c_1$ guarantees the positivity of $b_{22}(s_1)=q_2^{-1}U_2(s_1)$.
\end{remark}
We consider now the consistency condition \eqref{2tcss}. Combining \eqref{c1n} and \eqref{c2} gives
\begin{align*}
pc_1+\epsilon\kappa c_2=(c_1-\frac{2\epsilon}{n+1}s_1^{n+1})(p-\epsilon\kappa s_1)<0.
\end{align*}
Thus it follows from \eqref{2tcss} that
\begin{align}\label{c2n}
c_2<0.
\end{align}
\begin{remark}
The negativity of $c_2$ guarantees the positivity of $\Delta(s)$ on $[s_1,\infty)$. In fact, $\Delta(s_1)>0$ as shown above, and by \eqref{de}
\begin{eqnarray*}
\dot{\Delta}(s)=\frac{2p\kappa}{n+1}-\kappa^2c_2ns^{-n-1}>0,&&\forall s\ge s_1>0.
\end{eqnarray*}
\end{remark}
By \eqref{c2}, \eqref{c2n} is equivalent to
\begin{align}\label{c1p}
c_1>\frac{2\epsilon}{n+1}s_1^{n+1}-\frac{2p}{\kappa(n+1)}s_1^n.
\end{align}

\begin{lemma}\label{pa}
$\alpha>0$ on $(s_1,\infty)$.
\end{lemma}
\begin{proof}
By \eqref{aq}, we only need to show that $Q>0$ on $(s_1,\infty)$. Notice that
\begin{eqnarray*}
\dot{Q}(s)&=&-2\epsilon\frac{n+2}{n+1}s^{n+1}+\frac{2p}{\kappa}s^n+c_1\\
&>&-\frac{2\epsilon}{n+1}s^{n+1}+\frac{2p}{\kappa(n+1)}s^n+\frac{2\epsilon}{n+1}s_1^{n+1}-\frac{2p}{\kappa(n+1)}s_1^n\\
&=&-\frac{2\epsilon}{n+1}(s^{n+1}-s_1^{n+1})+\frac{2p}{\kappa(n+1)}(s^n-s_1^n)\\
&>&0,\;\;\;\;\forall s>s_1>0.
\end{eqnarray*}
The positivity of $Q(s)$ on $(s_1,\infty)$ then follows from the fact $Q(s_1)=0$.
\end{proof}
\begin{remark}
Now we take $I=(s_1,\infty)$. By Lemma \ref{pa}, there are no more collapses on $(s_1,\infty)$. So the metric $\widehat{g}$ (cf. \eqref{em}) is well-defined on $I$.
\end{remark}

We compute now the derivative of $\sqrt{b_{11}(t)}$ at $t_1$.
\begin{eqnarray*}
\frac{d}{dt}|_{t=t_1}\sqrt{b_{11}(t)}&=&\lim_{t\rightarrow t_1}\frac{\sqrt{b_{11}(t)}}{t-t_1}\\
&=&\lim_{t\rightarrow t_1}\frac{b_{11}'(t)}{2\sqrt{b_{11}(t)}}\\
&=&\lim_{s\rightarrow s_1}\frac{\dot{b_{11}}(s)\sqrt{\alpha(s)}}{2\sqrt{b_{11}(s)}}\\
&=&\frac{|q_2|\dot{\alpha}(s_1)}{2\sqrt{\Delta(s_1)}}.
\end{eqnarray*}
The smoothness conditions require that $|q_2|\dot{\alpha}(s_1)=2\sqrt{\Delta(s_1)}$, i.e.
\begin{eqnarray}
\label{2ttp}|q_2|=\frac{2\kappa^2s_1^{\frac{n+1}{2}}\sqrt{\frac{2\epsilon}{n+1}s_1^{n+1}-c_1}}{2s_1^n(p-\epsilon\kappa s_1)-\kappa(\frac{2\epsilon}{n+1}s_1^{n+1}-c_1)}\in\mathbb{Z}_+.
\end{eqnarray}
By \eqref{c1n} and \eqref{c1p}, we have
\begin{align*}
\frac{2\epsilon}{n+1}s_1^{n+1}-\frac{2p}{\kappa(n+1)}s_1^n<c_1<\frac{2\epsilon}{n+1}s_1^{n+1}.
\end{align*}
We may introduce a new parameter $0<\lambda<1$, and write
\begin{align*}
c_1=\frac{2\epsilon}{n+1}s_1^{n+1}-\frac{2p\lambda}{\kappa(n+1)}s_1^n.
\end{align*}
Thus (\ref{2ttp}) becomes
\begin{eqnarray}
\label{rf2ttp}|q_2|=\frac{\sqrt{\frac{2\lambda p\kappa^3s_1}{n+1}}}{p(1-\frac{\lambda}{n+1})-\epsilon\kappa s_1}\in\mathbb{Z}_+.
\end{eqnarray}
\begin{lemma}\label{kappa}
Given $s_1>0$ and $\lambda\in(0,1)$, there exists a unique $\kappa>0$ such that \eqref{rf2ttp} holds.
\end{lemma}

\begin{lemma}
$b_{ij}(t)$ and $c(t)$ are smooth and even functions in $t$ around $t_1$.
\end{lemma}
\begin{proof}
The smoothness of $b_{ij}(t)$ and $c(t)$ follows from the relation $\frac{d}{dt}=\sqrt{\alpha(s)}\frac{d}{ds}$ and the fact that $\alpha(s)$, $\Delta(s)$, $U_1(s)$, $U_2(s)$, and $\beta(s)$ are rational functions of $s$. The evenness property follows from the identity $\frac{d}{dt}|_{t=t_1}=\sqrt{\alpha(s)}\frac{d}{ds}|_{s=s_1}=0$ as $\alpha(s_1)=0$.
\end{proof}
Now we can extend the metric $\widehat{g}$ (cf. \eqref{em}) smoothly to the left endpoint $s_1$ of $I$ by choosing $\kappa$ as in Lemma \ref{kappa}. As noted before, this indicates the extended metric $g$ to be an Einstein metric.

\begin{proposition}\label{1}
We construct two two-parameter families of non-positive Einstein metrics with Einstein constant $\epsilon\le0$ on the associated solid torus bundle $T(P_q)$ of the form
\begin{eqnarray}\label{cneg}
g=\alpha(s)^{-1}ds^2+\sum_{i,j=1}^2b_{ij}(s)\theta^i\otimes\theta^j+\beta(s)\widehat{\pi}^\ast h,&s\in[s_1,\infty),
\end{eqnarray}
with free parameters $s_1>0$ and $\lambda\in(0,1)$. In details, $\kappa>0$ satisfies (\ref{rf2ttp}), and
\begin{align*}
\alpha=-\frac{2\epsilon}{n+1}s^2+\frac{2p}{\kappa(n+1)}s+c_1s^{1-n}+c_2s^{-n},&&\beta=\kappa s,&&\Delta=\frac{2p\kappa}{n+1}s+\kappa^2c_2s^{-n},
\end{align*}
\begin{align*}
&U_1=\frac{4(\lambda-1)pq_2\kappa^n(p-\epsilon\kappa s_1)s_1^n}{\psi(n+1)}(s_1^{n+1}s^{-n}-s),\\
&U_2=\frac{4(\lambda-1)pq_1\kappa^n(p-\epsilon\kappa s_1)s_1^n}{\psi(n+1)}(s-s_1^{n+1}s^{-n})+\frac{2p\kappa}{q_2(n+1)}(s+(\lambda-1)s_1^{n+1}s^{-n}),
\end{align*}
\begin{align*}
b_{11}=\frac{U_1^2+q_2^2\alpha}{\Delta},&&b_{12}=\frac{U_1U_2-q_1q_2\alpha}{\Delta},
&&b_{22}=\frac{U_2^2+q_1^2\alpha}{\Delta},
\end{align*}
\begin{align*}
c_1=\frac{2\epsilon}{n+1}s_1^{n+1}-\frac{2p\lambda}{\kappa(n+1)}s_1^n,&&
c_2=\frac{2p(\lambda-1)}{\kappa(n+1)}s_1^{n+1},
\end{align*}
\begin{align*}
\psi=\pm\sqrt{8(n+1)^{-1}p^2\kappa^{2n+1}\lambda(1-\lambda)s_1^{2n+1}(p-\epsilon\kappa s_1)}.
\end{align*}
\end{proposition}

\subsection{$\psi=0$}

In this case, \eqref{2tcs} reduces to $\dot{U}_1U_2-U_1\dot{U}_2=0$.
\begin{claim}\label{u1}
$U_1\equiv0$.
\end{claim}
\begin{proof}
If $U_1\ne0$, then we have
\begin{align*}
(\frac{U_2}{U_1})\dot{}=\frac{\dot{U}_2U_1-U_2\dot{U}_1}{U_1^2}=0.
\end{align*}
Thus $U_2=\mu U_1$ for some constant $\mu$. By Remark \ref{hhh}, we have $U_2(s_1)=\mu U_1(s_1)=0$. So $b_{22}(s_1)=0$ as $b_{12}(s_1)=0$ and $q_2\ne0$. This gives a contradiction.
\end{proof}
By Claim \ref{u1}, we have $\Delta=q_2U_2$, i.e., $U_2=q_2^{-1}\Delta$. So it follows from Lemma \ref{b} that
\begin{align*}
b_{11}=q_2^2\frac{\alpha}{\Delta},&&b_{12}=-q_1q_2\frac{\alpha}{\Delta},&&b_{22}=q_1^2\frac{\alpha}{\Delta}+q_2^{-2}\Delta.
\end{align*}
We turn now to check \eqref{2tv1} and \eqref{2tv2}. Notice that
\begin{align*}
\dot{b_{11}}=q_2^2\frac{\dot{\alpha}\Delta-\alpha\dot{\Delta}}{\Delta^2}.
\end{align*}
Substituting it into \eqref{2tv1} gives
\begin{align*}
q_2(\dot{\alpha}\Delta-\alpha\dot{\Delta})=v_1\Delta.
\end{align*}
By \eqref{als} and \eqref{de}, this amounts to
\begin{align*}
w_1=-p^{-1}q_2\kappa^n(npc_1+\epsilon\kappa c_2(n+1)),&&c_2(pc_1+\epsilon\kappa c_2)=0.
\end{align*}
On the other hand,
\begin{align*}
\dot{b_{12}}=-q_1q_2\frac{\dot{\alpha}\Delta-\alpha\dot{\Delta}}{\Delta^2}.
\end{align*}
Substituting it into \eqref{2tv2} gives
\begin{align*}
-q_1(\dot{\alpha}\Delta-\alpha\dot{\Delta})=v_2\Delta,
\end{align*}
which amounts to
\begin{align*}
w_2=p^{-1}q_1\kappa^n(npc_1+\epsilon\kappa c_2(n+1)).
\end{align*}
Clearly $\psi=q_1w_1+q_2w_2=0$ as supposed.

Regarding the consistency condition $c_2(pc_1+\epsilon\kappa c_2)=0$, we divide our discussion into two subcases $pc_1+\epsilon\kappa c_2=0$ and $c_2=0$.

\subsubsection{$pc_1+\epsilon\kappa c_2=0$}

In this case, $c_1=-p^{-1}\epsilon\kappa c_2$. Thus
\begin{align*}
&\alpha=-\frac{2\epsilon}{n+1}s^2+\frac{2p}{\kappa(n+1)}s-p^{-1}\epsilon\kappa c_2s^{1-n}+c_2s^{-n}
=(\kappa^{-1}s^{-n}-p^{-1}\epsilon s^{1-n})(\frac{2p}{n+1}s^{n+1}+\kappa c_2),\\
&\Delta=\frac{2p\kappa}{n+1}s+\kappa^2c_2s^{-n}=\kappa s^{-n}(\frac{2p}{n+1}s^{n+1}+\kappa c_2).
\end{align*}
We have
\begin{eqnarray*}
b_{11}=q_2^2\frac{\alpha}{\Delta}=q_2^2\kappa^{-2}(1-p^{-1}\epsilon\kappa s)>0,&&\forall s\ge s_1\ge0.
\end{eqnarray*}
In particular, $b_{11}(s_1)>0$, which contradicts the smooth collapse conditions.

\subsubsection{$c_2=0$}

In this case, we have
\begin{eqnarray*}
\alpha=-\frac{2\epsilon}{n+1}s^2+\frac{2p}{\kappa(n+1)}s+c_1s^{1-n},&&\Delta=\frac{2p\kappa}{n+1}s.
\end{eqnarray*}
\begin{claim}\label{s1}
$s_1>0$.
\end{claim}
\begin{proof}
If $s_1=0$, then $c_1$ has to be zero since $\alpha(s_1)=0$. Thus
\begin{eqnarray*}
b_{11}=q_2^2\frac{\alpha}{\Delta}=q_2^2\kappa^{-2}(1-p^{-1}\epsilon\kappa s)>0,&&\forall s\ge0.
\end{eqnarray*}
As shown above, this contradicts the smooth collapse conditions.
\end{proof}
An immediate consequence of Claim \ref{s1} is that
\begin{corollary}
$\Delta>0$ on $[s_1,\infty)$.
\end{corollary}
On the other hand, the requirement $\alpha(s_1)=0$ is equivalent to
\begin{align}\label{c11}
c_1=\frac{2\epsilon}{n+1}s_1^{n+1}-\frac{2p}{\kappa(n+1)}s_1^n.
\end{align}
Notice that $c_1<0$.
\begin{lemma}\label{a}
$\alpha>0$ on $(s_1,\infty)$.
\end{lemma}
\begin{proof}
By \eqref{c11}, we have
\begin{eqnarray*}
\dot{\alpha}=-\frac{4\epsilon}{n+1}s+\frac{2p}{\kappa(n+1)}+c_1(1-n)s^{-n}>0,&&\forall s\ge s_1>0.
\end{eqnarray*}
Now the claim follows from the fact $\alpha(s_1)=0$.
\end{proof}

Lemma \ref{a} implies that there are no more collapses after $s_1$. So we can take $I=(s_1,\infty)$. The metric $\widehat{g}$ (cf. \eqref{em}) is well-defined on $I$.

It follows from $\alpha(s_1)=0$ and $\Delta(s_1)>0$ that $b_{11}(s_1)=b_{12}(s_1)=0$ and $b_{22}(s_1)>0$. Moreover,
\begin{align*}
\frac{d}{dt}|_{t=t_1}\sqrt{b_{11}(t)}=\frac{|q_2|\dot{\alpha}(s_1)}{2\sqrt{\Delta(s_1)}}=\frac{|q_2|(-\epsilon\kappa s_1+\frac{np}{n+1})}{\sqrt{\frac{2p\kappa^3}{n+1}s_1}}=1,
\end{align*}
as required by the smooth collapse conditions, is equivalent to
\begin{eqnarray}\label{2tk0}
|q_2|=\frac{\sqrt{\frac{2p\kappa^3}{n+1}s_1}}{(-\epsilon\kappa s_1+\frac{np}{n+1})}\in\mathbb{Z}_+.
\end{eqnarray}
\begin{proposition}\label{ks}
Given $s_1>0$, there exists a unique $\kappa>0$ such that \eqref{2tk0} holds.
\end{proposition}
As before, it is easy to see that $b_{ij}(t)$ and $c(t)$ are smooth and even in $t$ around $t_1$. So by choosing $\kappa$ as in Proposition \ref{ks}, we can extend the metric $\widehat{g}$ (cf. \eqref{em}) smoothly to the left endpoint $s_1$ of $I$. The extended metric $g$ is an Einstein metric.

\begin{proposition}\label{2} We construct a one-parameter family of non-positive Einstein metrics with Einstein constant $\epsilon\le0$ on the associated solid torus bundle $T(P_q)$ of the form
\begin{eqnarray*}
g=\alpha(s)^{-1}ds^2+\sum_{i,j=1}^2b_{ij}(s)\theta^i\otimes\theta^j+\beta(s)\widehat{\pi}^\ast h,&&s\in[s_1,\infty),
\end{eqnarray*}
with free parameter $s_1>0$. In details, $\kappa>0$ satisfies \eqref{2tk0}, and
\begin{align*}
\alpha=-\frac{2\epsilon}{n+1}s^2+\frac{2p}{\kappa(n+1)}s+(\frac{2\epsilon}{n+1}s_1^{n+1}-\frac{2p}{\kappa(n+1)}s_1^n)s^{1-n},&&\beta=\kappa s,&&\Delta=\frac{2p\kappa}{n+1}s,
\end{align*}
\begin{align*}
b_{11}=q_2^2\frac{\alpha}{\Delta},&&b_{12}=-q_1q_2\frac{\alpha}{\Delta},&&b_{22}=q_1^2\frac{\alpha}{\Delta}+q_2^{-2}\Delta.
\end{align*}
\end{proposition}

\begin{remark}
The boundary one-parameter family with $\lambda=1$ in Proposition \ref{1} corresponds to the one in the preceding proposition. So we may combine Propositions \ref{1} and \ref{2} by allowing $\lambda=1$ in the former one.
\end{remark}

\subsection{CCE manifolds}

In this subsection, we specialize to the negative Einstein metrics in Proposition \ref{1} (including the limiting case $\lambda=1$). We will show that they are CCE metrics, and thus complete the proof of Theorem \ref{neg}. Without loss of generality, we assume $\epsilon=-(2n+2)$.
\begin{lemma}
The Einstein metric $g$ given by \eqref{cneg} is complete.
\end{lemma}
\begin{proof}
The geodesic distance $\int_{s_1}^s\alpha(\tau)^{-\frac{1}{2}}d\tau$ w.r.t. $g$ is unbounded since $\alpha(s)\sim s^2$ as $s\rightarrow\infty$.
\end{proof}
We define now a conformal metric $\overline{g}$ on $T(P_q)$ given by
\begin{align}\label{cm}
\overline{g}=s^{-1}g=(s\alpha(s))^{-1}ds^2+\sum_{i,j=1}^2s^{-1}b_{ij}(s)\theta^i\otimes\theta^j+\kappa\widehat{\pi}^\ast h.
\end{align}
\begin{lemma}
$\overline{g}$ extends to a compact metric on $\overline{T(P_q)}=T(P_q)\cup\partial T(P_q)$.
\end{lemma}
\begin{proof}
The geodesic distance $\int_{s_1}^s(\tau\alpha(\tau))^{-\frac{1}{2}}d\tau$ w.r.t. $\overline{g}$ is bounded since $\alpha(s)\sim s^2$ as $s\rightarrow\infty$. Notice that $U_i(s)$ and $\Delta(s)$ grow like $s$ as $s\rightarrow\infty$. So by Proposition \ref{1}, $b_{ij}(s)\sim s$ as $s\rightarrow\infty$. Moreover, $\lim_{s\rightarrow\infty}s^{-2}\alpha(s)=4>0$. Thus $\overline{g}|_{\partial T(P_q)=\{s=\infty\}}$ is well-defined and non-degenerate.
\end{proof}
\begin{corollary}
$g$ is a CCE metric.
\end{corollary}
\begin{remark}
Note that $\rho(s)=s^{-\frac{1}{2}}$ is a defining function for $\partial T(P_q)$. In fact, it is clear that $\rho>0$ in $T(P_q)=\{s_1\le s<\infty\}$, and $\rho=0$ on $\partial T(P_q)=\{s=\infty\}$. Moreover,
\begin{align*}
\frac{\partial}{\partial\rho}=\frac{ds}{d\rho}\frac{\partial}{\partial s}=-2\rho^{-3}\frac{\partial}{\partial s}=-2s^\frac{3}{2}\frac{\partial}{\partial s}.
\end{align*}
Thus
\begin{align*}
\lim_{s\rightarrow\infty}\overline{g}(\frac{\partial}{\partial\rho},\frac{\partial}{\partial\rho})=\lim_{s\rightarrow\infty}4s^2\alpha(s)^{-1}=1>0.
\end{align*}
In particular, $d\rho|_{\partial T(P_q)}\ne0$.
\end{remark}
\begin{remark}
The event horizon $\{s=s_1\}$ is the nontrivial principal circle bundle over $V$ with Euler class $q_2a$. Topologically, it is the deformation retraction of $T(P_q)$,
\end{remark}

\subsection{$Q$-curvatures}\label{ci}

In this subsection, we will give a proof of Theorem \ref{zq}. The proof consists of two steps. First, we show the vanishing of total $Q$-curvature for every conformal infinity associated with the CCE manifolds in Theorem \ref{neg}. Second, we show that every conformal infinity contains a representative with constant $Q$-curvature, which thus must have constant zero $Q$-curvature.

\subsubsection{Geodesic defining functions}\label{gdf}

Recall that the CCE metrics in Theorem \ref{neg} are given by \eqref{cneg}. We see that $\rho=s^{-\frac{1}{2}}$ is not a geodesic defining function. In order to find a geodesic defining function $\sigma$, we are led to a first-order ordinary differential equation
\begin{align}\label{fd}
\frac{d\sigma}{\sigma}=-\frac{ds}{\sqrt{\alpha(s)}}.
\end{align}

\begin{proposition}\label{ggdf}
Let
\begin{eqnarray*}
\zeta(s)=\int_{s_1}^s(\frac{1}{2\tau}-\frac{1}{\sqrt{\alpha(\tau)}})d\tau,&&s\ge s_1>0.
\end{eqnarray*}
Then $\sigma=s^{-\frac{1}{2}}\exp\zeta$ satisfies \eqref{fd}, and is a geodesic defining function.
\end{proposition}
\begin{proof}
At $s_1$, $\alpha(s_1)=0$ and $\alpha'(s_1)=4(n+1)s_1+2p\kappa^{-1}(1-\frac{\lambda}{n+1})>0$. Thus $\alpha(s)=O(|s-s_1|)$. We see that $\zeta$ is well-defined on $[s_1,\infty)$. In particular, $\zeta(s_1)=0$. On the other hand,
\begin{align*}
\lim_{s\rightarrow\infty}\frac{\alpha(s)-4s^2}{s}=\frac{2p}{\kappa(n+1)}>0,
\end{align*}
i.e., $\alpha(s)-4s^2\sim s$ as $s\rightarrow\infty$. Thus
\begin{align*}
\frac{1}{2\tau}-\frac{1}{\sqrt{\alpha(\tau)}}=\frac{\alpha(\tau)-4\tau^2}{2\tau\sqrt{\alpha(\tau)}(\sqrt{\alpha(\tau)}+2\tau)}=O(\frac{1}{\tau^2}).
\end{align*}
We see that $\zeta_{\infty}=\lim_{s\rightarrow\infty}\zeta(s)$ is finite. By definition, $\lim_{s\rightarrow\infty}\sigma=0$, i.e., $\sigma=0$ on $\partial T(P_q)$.

Now differentiating the formula of $\zeta$ gives
\begin{eqnarray}\label{zet}
\dot{\zeta}=\frac{1}{2s}-\frac{1}{\sqrt{\alpha(s)}},&&s>s_1,
\end{eqnarray}
which shows that $\zeta$ is smooth on $(s_1,\infty)$.

Finally, it is not hard to check that $\sigma=s^{-\frac{1}{2}}\exp\zeta$ satisfies \eqref{fd}. Thus $\sigma$ is a geodesic defining function. This completes the proof of Proposition \ref{ggdf}.
\end{proof}

\begin{remark}
$\zeta$ is not smooth at $s_1$. To see this, we return to $t$-coordinate. By \eqref{zet}, we have
\begin{align*}
\frac{d\zeta}{dt}=\frac{d\zeta}{ds}\frac{ds}{dt}=\frac{\sqrt{\alpha}}{2s}-1.
\end{align*}
In particular, $\dot{\zeta}(t_1)=-1\ne0$. Thus $\zeta$ is not even in $t$ around $t_1$, and hence not smooth at $t_1$.
\end{remark}

By \eqref{fd},
\begin{eqnarray*}
\frac{d\sigma}{ds}=-\frac{\sigma}{\sqrt{\alpha(s)}}<0,&&s>s_1.
\end{eqnarray*}
Hence the range of $\sigma$ is $[0,s_1^{-\frac{1}{2}}]$ since $\sigma(s_1)=s_1^{-\frac{1}{2}}$ as $\zeta(s_1)=0$, and $\lim_{s\rightarrow\infty}\sigma=0$. By the inverse function theorem, $s$ can be written as a function of $\sigma$ on $[0,s_1^{-\frac{1}{2}}]$, say $s=\sigma^{-2}\xi$, where $\xi=\xi(\sigma)$ is the unique smooth solution of the ordinary differential equation
\begin{eqnarray*}
\frac{d\xi}{d\sigma}&=&2\sigma^{-1}\xi(1-\sqrt{\frac{\alpha(s)}{4s^2}})\\
&=&-\frac{\frac{2p}{\kappa(n+1)}\sigma+c_1\sigma^{2n+1}\xi^{-n}+c_2\sigma^{2n+3}\xi^{-n-1}}{2+\sqrt{4+\frac{2p}{\kappa(n+1)}\sigma^2\xi^{-1}+c_1\sigma^{2n+2}\xi^{-n-1}+c_2\sigma^{2n+4}\xi^{-n-2}}},
\end{eqnarray*}
subject to the initial condition $\xi(0)=\exp(2\zeta_{\infty})$.

\subsubsection{Asymptotic volumes}\label{aav}

In terms of the geodesic defining function $\sigma$, we can rewrite our Einstein metric $g$ as
\begin{align}\label{gm}
g=\sigma^{-2}d\sigma^2+\sum_{i,j=1}^2E_{ij}(\sigma)\theta^i\otimes\theta^j+F(\sigma)\widehat{\pi}^\ast h
\end{align}
with $E_{ij}(\sigma)=b_{ij}(s)$ and $F(\sigma)=\beta(s)$. By \eqref{fd} and \eqref{gm}, the asymptotic volume is
\begin{eqnarray*}
\mbox{Vol}_g(\{\sigma>\delta\})&=&C\int_\delta^{\sigma(s_1)}\sqrt{\det(E_{ij})}F^n\frac{d\sigma}{\sigma}\\
&=&C\int_{s(\delta)}^{s_1}\sqrt{\det(b_{ij})}\beta^n(-\frac{ds}{\sqrt{\alpha}})\\
&=&C\kappa^n\int_{s_1}^{s(\delta)}s^nds\\
&=&\frac{C\kappa^n}{n+1}(s(\delta)^{n+1}-s_1^{n+1})\\
&=&\frac{C\kappa^n}{n+1}(\delta^{-(2n+2)}\xi(\delta)^{n+1}-s_1^{n+1}),
\end{eqnarray*}
where $C=4\pi^2\mbox{Vol}_h(V)$. Since $\xi$ is smooth around $0$, no $\log\delta$ term will appear in the asymptotic expansion of $\mbox{Vol}_g(\{r>\delta\})$. Thus the conformal infinity has vanishing total $Q$-curvature.

\subsubsection{Constant $Q$-curvature metrics}\label{cq}

Now it remains to show the existence of representatives with constant $Q$-curvature in every conformal infinity. Let us consider the restriction $g_b$ of $\sigma^2g$ to the boundary $\partial T(P_q)$. To show that $g_b$ has constant $Q$-curvature, by \cite[Theorem 3.1]{FefGra02}, we need to find the unique solution $\widetilde{U}$ mod $O(\sigma^{2n+2})$ of the Poisson equation
\begin{align*}
\triangle_g\widetilde{U}=2n+2+O(\sigma^{2n+3}\log\sigma)
\end{align*}
of the form
\begin{align*}
\widetilde{U}=\log\sigma +A+B\sigma^{2n+2}\log\sigma+O(\sigma^{2n+2}),
\end{align*}
with $A,B\in C^\infty(T(P_q))$ and $A|_{\partial T(P_q)}=0$. It has been shown that $B|_{\partial T(P_q)}$ is a constant multiple of the $Q$-curvature of $g_b$ \cite[Theorem 4.1]{FefGra02}.
\begin{lemma}
For any $C^2$ function $G(x)$ of variable $x$, $\triangle_g(G(\sigma))$ depends on $\sigma$ only.
\end{lemma}
\begin{proof}
Let $\iota=\log\sigma$, i.e., $\sigma=\exp\iota$. Thus $G(\sigma)=G(\log\iota)$. We may rewrite our Einstein metric as
\begin{align*}
g=d\iota^2+\sum_{i,j=1}^2\overline{E}_{ij}(\iota)\theta^i\otimes\theta^j+\overline{F}(\iota)\pi^\ast h,
\end{align*}
with $\overline{E}_{ij}(\iota)=E_{ij}(\sigma)$ and $\overline{F}(\iota)=F(\sigma)$. Then it is not hard to see that $\triangle_g(G(\log\iota))$ is a function of $\iota$ only (cf. Lemma \ref{sh}). In other words, $\triangle_g(G(\sigma))$ is a function of $\sigma$ only.
\end{proof}
Following closely the proof of \cite[Theorem 3.1]{FefGra02}, we first construct the function $A$ mod $O(\sigma^{2n+2})$ of the form $\sum_{i=2}^{2n+1}c_i\sigma^i$ with constants $c_i$ such that
\begin{align*}
\triangle_gA=2n+2-\triangle_g\log\sigma+\sigma^{2n+2}\widetilde{W},
\end{align*}
where the remainder $\widetilde{W}$ turns out to depend on $\sigma$ only. Then we can formally determine $B=B(\sigma)$ subject to the boundary condition $B|_{\partial T(P_q)}=-(2n+2)^{-1}\widetilde{W}(0)$, from which follows the constancy of the $Q$-curvature of $g_b$. As shown in \S\ref{aav}, $g_b$ has vanishing total $Q$-curvature. Thus $g_b$ actually has vanishing $Q$-curvature. This completes the proof of Theorem \ref{zq}.

\begin{remark}
The $Q$-curvature in dimension $4$ is given by an explicit formula
\begin{align*}
Q=\frac{1}{6}(\triangle R+R^2-3|\mbox{Ric}|^2).
\end{align*}
So it is not hard to verify that for every $5$-dimensional CCE manifold in Theorem \ref{neg}, the restriction of $\overline{g}$ (cf. \eqref{cm}) to $\partial T(P_q)=\{s=\infty\}$ has zero $Q$-curvature (cf. Lemma \ref{ip} and Proposition \ref{1}).
\end{remark}

\subsection{Ricci-flat manifolds}

In this subsection, we specialize to the Ricci-flat metrics in Proposition \ref{1} (including again the limiting case $\lambda=1$). We will show that they are complete, and have slow volume growth and quadratic curvature decay.

\subsubsection{Completeness}

The geodesic distance w.r.t. $g$ (cf. \eqref{cneg}) is given by
\begin{align*}
t=\int\alpha(s)^{-\frac{1}{2}}ds.
\end{align*}
By Proposition \ref{1}, $\alpha(s)\sim s$ as $s\rightarrow\infty$. So $t\sim\sqrt{s}$ is unbounded. Thus $g$ is complete.
\begin{remark}
More precisely, we have $t\sim\sqrt{2p^{-1}\kappa(n+1)s}$. By Proposition \ref{1}, $b_{ij}(s)$ and $\beta(s)$ grow like $s$ as $s\rightarrow\infty$. In other words, $b_{ij}(t)$ and $c(t)^2$ grow like $t^2$. So $g$ is asymptotic to a cone metric
\begin{align*}
g_c=dt^2+t^2(\sum_{i,j=1}^2\widetilde{b}_{ij}\theta^i\otimes\theta^j+\frac{p}{2(n+1)}\widehat{\pi}^\ast h),
\end{align*}
where $\widetilde{b}_{ij}=\frac{1}{4\kappa^2}\widetilde{U}_i\widetilde{U}_j$, and
\begin{align*}
\widetilde{U}_1=-\frac{4(\lambda-1)p^2q_2\kappa^ns_1^n}{\psi(n+1)},&&\widetilde{U}_2=\frac{4(\lambda-1)p^2q_1\kappa^ns_1^n}{\psi(n+1)}+\frac{2p\kappa}{q_2(n+1)}.
\end{align*}
However, $g_c$ is degenerate in the sense that $\det(\widetilde{b}_{ij})=0$.
\end{remark}

\subsubsection{Volume growth}

The volume element w.r.t. $g$ is given by
\begin{align*}
\mbox{Vol}_g=\sqrt{\frac{1}{\alpha}\cdot\alpha\cdot\beta^{2n}}\cdot\mbox{Vol}_h\sim s^n.
\end{align*}
So the volume function grows like $s^{n+1}$, i.e., $t^{2n+2}$ since $s\sim t^2$. Thus $g$ has slower-than-Euclidean volume growth.

\subsubsection{Curvature decay}

Using formulas from \S\ref{fpt}, one can check that the curvature tensor of $g$ decays as fast as $t^{-2}$ when $t\rightarrow\infty$, and hence so do the sectional curvatures of $g$.

\section{Positive Einstein manifolds}\label{pe}

In this section, we prove Theorem \ref{pos} by constructing complete positive Einstein metrics from local solutions given in Proposition \ref{ls}. From now on, we fix the Einstein constant $\epsilon=2n+2$.

\subsection{Smooth collapse}

In order to obtain closed manifolds from $\widehat{M}=I\times P_q$, we must add to it two principal circle bundles, one at each endpoint of $I$. In other words, we shall collapse the two standard circles in each $2$-torus fiber of $P_q$, one at each endpoint of $I$. More precisely, we collapse the one corresponding to characteristic class $q_1a$ at the left endpoint $s_1$, while the one corresponding to characteristic class $q_2a$ at the right endpoint $s_2$. The smoothness conditions on $b_{ij}(t)$'s are
\begin{enumerate}
  \item $b_{ij}(t)$'s are smooth and even in $t$ around $t_k=t(s_k)$, $k=1,2$;
  \item $b_{11}(t_1)=b_{12}(t_1)=0$, $b_{22}(t_1)>0$, and $\frac{d}{dt}|_{t=t_1}\sqrt{b_{11}(t)}=1$;
  \item $b_{12}(t_2)=b_{22}(t_2)=0$, $b_{11}(t_2)>0$, and $\frac{d}{dt}|_{t=t_2}\sqrt{b_{22}(t)}=-1$.
\end{enumerate}

\begin{remark}
Recall that $U_1=q_1b_{11}+q_2b_{12}$, $U_2=q_1b_{12}+q_2b_{22}$, and $\alpha=\det B$ given by \eqref{als}
\begin{align}\label{alp}
\alpha(s)=-4s^2+\frac{2p}{\kappa(n+1)}s+c_1s^{1-n}+c_2s^{-n}.
\end{align}
Conditions $(2)$ and $(3)$ imply $\alpha(s_i)=U_i(s_i)=0$, $i=1,2$.
\end{remark}

In \S\ref{e}, local solutions fall into two categories, $\psi\ne0$ and $\psi=0$. In the latter case, we can apply the argument for Claim \ref{u1} to the current situation, and conclude that $U_i\equiv0$, $i=1,2$, since $|q_1|>|q_2|>0$ by assumption. It follows that $\Omega\equiv0$, which gives a contradiction. Thus we cannot construct smooth positive Einstein metrics when $\psi=0$. Henceforth, we focus on the case $\psi\ne0$.

We proceed to study the smooth collapse at $s_1\ge0$. As shown at the start of \S\ref{pn}, the assumption $\psi\ne0$ forces $s_1>0$. By \eqref{alp}, $\alpha(s_1)=0$ is equivalent to
\begin{align}\label{c21}
c_2=4s_1^{n+2}-\frac{2p}{\kappa(n+1)}s_1^{n+1}-c_1s_1.
\end{align}
Substituting \eqref{c21} into $U_1(s_1)=0$ (cf. \eqref{2tue}) yields
\begin{align*}
(4s_1^{n+1}-c_1)(w_1-q_2\kappa^{n-1}(\kappa c_1+2ps_1^n-(4n+4)\kappa s_1^{n+1}))=0.
\end{align*}
\begin{claim}
$c_1\ne4s_1^{n+1}$.
\end{claim}
\begin{proof}
If $c_1=4s_1^{n+1}$, then \eqref{c21} reduces to $c_2=-\frac{2p}{\kappa(n+1)}s_1^{n+1}$. So $pc_1+(2n+2)\kappa c_2=0$. By \eqref{2tcss}, $\psi^2=0$. This gives a contradiction.
\end{proof}
Therefore,
\begin{align}\label{w11}
w_1=q_2\kappa^{n-1}(\kappa c_1+2ps_1^n-(4n+4)\kappa s_1^{n+1}).
\end{align}
\begin{remark}
\eqref{w11} necessitates $q_2\ne0$. If $q_2=0$, then $w_1=0$, and hence $\psi=q_1w_1+q_2w_2=0$. This gives a contradiction.
\end{remark}

By \eqref{de} and \eqref{c21},
\begin{align*}
\Delta(s_1)=\frac{2p\kappa}{n+1}s_1+\kappa^2s_1^{-n}(4s_1^{n+2}-\frac{2p}{\kappa(n+1)}s_1^{n+1}-c_1s_1)=\kappa^2s_1^{1-n}(4s_1^{n+1}-c_1)\ge0.
\end{align*}
So $c_1\le4s_1^{n+1}$. On the other hand, $pc_1+(2n+2)\kappa c_2=(p-(2n+2)\kappa s_1)(c_1-4s_1^{n+1})$. Thus \eqref{2tcss} becomes
\begin{align}\label{ps}
\psi^2=(2n+2)\kappa^{2n+3}c_2(p-(2n+2)\kappa s_1)(c_1-4s_1^{n+1})>0.
\end{align}
In particular, $c_1-4s_1^{n+1}\ne0$. So it must be
\begin{align}\label{c1l}
c_1<4s_1^{n+1}.
\end{align}
Thus $\Delta(s_1)>0$. It follows that $U_2(s_1)=q_2^{-1}\Delta(s_1)\ne0$. By Lemma \ref{b}, we see that
\begin{proposition}
\eqref{c21}, \eqref{w11}, and \eqref{c1l} ensure $b_{11}(t_1)=b_{12}(t_1)=0$ and $b_{12}(t_1)>0$.
\end{proposition}
\begin{remark}
By \eqref{c1l}, \eqref{ps} amounts to
\begin{align}\label{c21l}
c_2(p-(2n+2)\kappa s_1)<0.
\end{align}
\end{remark}

We study now the smooth collapse at $s_2(>s_1>0)$. It follows from $\alpha(s_2)=0$ that
\begin{align}\label{c22}
c_2=4s_2^{n+2}-\frac{2p}{\kappa(n+1)}s_2^{n+1}-c_1s_2.
\end{align}
Substituting \eqref{c22} into $U_2(s_2)=0$ (cf. \eqref{2tve}) yields
\begin{align*}
(4s_2^{n+1}-c_1)(w_2+q_1\kappa^{n-1}(\kappa c_1+2ps_2^n-(4n+4)\kappa s_2^{n+1}))=0.
\end{align*}
By \eqref{c1l}, we have $c_1<4s_2^{n+1}$. Therefore,
\begin{align}\label{w22}
w_2=-q_1\kappa^{n-1}(\kappa c_1+2ps_2^n-(4n+4)\kappa s_2^{n+1}).
\end{align}
\begin{remark}
\eqref{w22} necessitates $q_1\ne0$. If $q_1=0$, then $w_2=0$, and hence $\psi=q_1w_1+q_2w_2=0$. This gives a contradiction.
\end{remark}

By \eqref{de} and \eqref{c22},
\begin{align*}
\Delta(s_2)=\frac{2p\kappa}{n+1}s_2+\kappa^2s_2^{-n}(4s_2^{n+2}-\frac{2p}{\kappa(n+1)}s_2^{n+1}-c_1s_2)=\kappa^2s_2^{1-n}(4s_2^{n+1}-c_1)>0.
\end{align*}
Thus $U_1(s_2)=q_1^{-1}\Delta(s_2)\ne0$. By Lemma \ref{b}, we see that
\begin{proposition}
\eqref{c22}, \eqref{w22}, and \eqref{c1l} ensure $b_{12}(t_2)=b_{22}(t_2)=0$ and $b_{11}(t_2)>0$.
\end{proposition}
By \eqref{c22}, we have $pc_1+(2n+2)\kappa c_2=(p-(2n+2)\kappa s_2)(c_1-4s_2^{n+1})$. So \eqref{2tcss} becomes
\begin{align*}
\psi^2=(2n+2)\kappa^{2n+3}c_2(p-(2n+2)\kappa s_2)(c_1-4s_2^{n+1})>0.
\end{align*}
Thus
\begin{align}\label{c22l}
c_2(p-(2n+2)\kappa s_2)<0.
\end{align}

\begin{lemma}
$p>(2n+2)\kappa s_1$.
\end{lemma}
\begin{proof}
Let
\begin{align}\label{fz}
F(z)=4z^{n+2}-\frac{2p}{\kappa(n+1)}z^{n+1}-c_1z,&s_1\le z\le s_2.
\end{align}
By \eqref{c21} and \eqref{c22}, $F(s_1)=F(s_2)$. Thus Rolle's theorem asserts a $z_0\in(s_1,s_2)$ such that $F'(z_0)=0$. However,
\begin{eqnarray*}
F'(z)&=&4(n+2)z^{n+1}-2p\kappa^{-1}z^n-c_1\\
&>&4(n+2)z^{n+1}-2p\kappa^{-1}z^n-4s_1^{n+1}\\
&=&4(z^{n+1}-s_1^{n+1})+2\kappa^{-1}z^n((2n+2)\kappa z-p)\\
&>&2\kappa^{-1}z^n((2n+2)\kappa s_1-p),\;\;\;\;s_1<z<s_2.
\end{eqnarray*}
If $p\le(2n+2)\kappa s_1$, then $F'(z)>0$ on $(s_1,s_2)$. This gives a contradiction. Thus $p>(2n+2)\kappa s_1$.
\end{proof}
By \eqref{c21l} and \eqref{c22l}, $c_2<0$ and $p>(2n+2)\kappa s_2$. Therefore,
\begin{align}\label{s12}
0<s_1<s_2<\frac{p}{\kappa(2n+2)}.
\end{align}
On the other hand, equating \eqref{c21} with \eqref{c22} gives
\begin{align}\label{c1}
c_1=4\frac{s_2^{n+2}-s_1^{n+2}}{s_2-s_1}-\frac{2p}{\kappa(n+1)}\frac{s_2^{n+1}-s_1^{n+1}}{s_2-s_1}.
\end{align}
Substituting \eqref{c1} into \eqref{c21} gives
\begin{align}\label{c2f}
c_2=\frac{s_1s_2}{s_1-s_2}(4(s_2^{n+1}-s_1^{n+1})-\frac{2p}{\kappa(n+1)}(s_2^n-s_1^n)).
\end{align}
Since $c_2<0$, \eqref{c21} and \eqref{c22} imply
\begin{eqnarray}\label{c1g}
c_1>4s_k^{n+1}-\frac{2p}{\kappa(n+1)}s_k^n,&&k=1,2.
\end{eqnarray}
\begin{lemma}
Suppose $c_1$ is given by \eqref{c1}. Then $c_1<4s_1^{n+1}$ provided $s_2<\frac{p}{\kappa(2n+2)}$. Furthermore, \eqref{c1g} is equivalent to
\begin{align}\label{do}
s_2^{n+1}-\frac{p}{\kappa(2n+2)}s_2^n>s_1^{n+1}-\frac{p}{\kappa(2n+2)}s_1^n.
\end{align}
\end{lemma}
\begin{proof}
By \eqref{c1},
\begin{align*}
c_1-4s_1^{n+1}=4\frac{s_2^{n+1}-s_1^{n+1}}{s_2-s_1}(s_2-\frac{p}{\kappa(2n+2)}).
\end{align*}
So $c_1-4s_1^{n+1}<0$ if $s_2-\frac{p}{\kappa(2n+2)}<0$. Furthermore,
\begin{align*}
&&c_1-4s_1^{n+1}+\frac{2p}{\kappa(n+1)}s_1^n=\frac{4s_2}{s_2-s_1}(s_2^{n+1}-s_1^{n+1}-\frac{p}{\kappa(2n+2)}(s_2^n-s_1^n)),\\
&&c_1-4s_2^{n+1}+\frac{2p}{\kappa(n+1)}s_2^n=\frac{4s_1}{s_2-s_1}(s_2^{n+1}-s_1^{n+1}-\frac{p}{\kappa(2n+2)}(s_2^n-s_1^n)).
\end{align*}
They are positive iff $s_2^{n+1}-s_1^{n+1}-\frac{p}{\kappa(2n+2)}(s_2^n-s_1^n)>0$.
\end{proof}

\begin{lemma}\label{f}
$F(z)<c_2<0$ on $(s_1,s_2)$.
\end{lemma}
\begin{proof}
By \eqref{fz},
\begin{align*}
F''(z)=4(n+1)(n+2)z^{n-1}(z-\frac{np}{2\kappa(n+1)(n+2)}).
\end{align*}
Thus
\begin{align*}
F''(z)\left\{\begin{array}{cl}
<0,&z\in(0,\frac{np}{2\kappa(n+1)(n+2)});\\
>0,&z\in(\frac{np}{2\kappa(n+1)(n+2)},\infty),
\end{array}\right.
\end{align*}
Since $F(0)=0>c_2=F(s_1)=F(s_2)$, it is not hard to see that $F(z)<c_2$ on $(s_1,s_2)$.
\end{proof}

\begin{corollary}\label{ap}
$\alpha(s)>0$ on $(s_1,s_2)$.
\end{corollary}
\begin{proof}
Assume $\alpha(s_0)\le0$ for some $s_0\in(s_1,s_2)$. Then $c_2\le4s_0^{n+2}-\frac{2p}{\kappa(n+1)}s_0^{n+1}-c_1s_0=F(s_0)$. This contradicts Lemma \ref{f}.
\end{proof}

\begin{lemma}\label{del}
$\Delta(s)>0$ on $(s_1,s_2)$.
\end{lemma}
\begin{proof}
By \eqref{de},
\begin{eqnarray*}
\dot{\Delta}=\frac{2p\kappa}{n+1}-n\kappa^2c_2s^{-n-1}>0,&&s>0.
\end{eqnarray*}
Thus $\Delta(s)>\Delta(s_1)>0$ when $s>s_1$.
\end{proof}

\begin{remark}
Lemmas \ref{b}, \ref{del}, and Corollary \ref{ap} imply that $B(s)=(b_{ij}(s))$ is well-defined and positive-definite on $(s_1,s_2)$.
\end{remark}

We study now remaining smoothness conditions at $s_i$, $i=1,2$. As shown before,
\begin{align*}
\frac{d}{dt}|_{t=t_1}\sqrt{b_{11}(t)}=\frac{|q_2|\dot{\alpha}(s_1)}{2\sqrt{\Delta(s_1)}}.
\end{align*}
So $\frac{d}{dt}|_{t=t_1}\sqrt{b_{11}(t)}=1$ iff $q_2^2=4\Delta(s_1)\dot{\alpha}(s_1)^{-2}$. However,
\begin{eqnarray*}
\dot{\alpha}(s_1)&=&s_1^{-n}(-4(n+2)s_1^{n+1}+2p\kappa^{-1}s_1^n+c_1)\\
&=&2\kappa^{-1}(p-(2n+2)\kappa s_1-\frac{(p-(2n+2)\kappa s_2)(s_2^{n+1}-s_1^{n+1})}{(n+1)s_1^n(s_2-s_1)}),
\end{eqnarray*}
and
\begin{align*}
\Delta(s_1)=4\kappa^2s_1^{1-n}(\frac{p}{\kappa(2n+2)}-s_2)\frac{s_2^{n+1}-s_1^{n+1}}{s_2-s_1}.
\end{align*}
Thus
\begin{align}\label{q2}
q_2^2=\frac{2\kappa^3(p-(2n+2)\kappa s_2)(s_2^{n+1}-s_1^{n+1})}{(n+1)s_1^{n-1}(s_2-s_1)(p-(2n+2)\kappa s_1-\frac{(p-(2n+2)\kappa s_2)(s_2^{n+1}-s_1^{n+1})}{(n+1)s_1^n(s_2-s_1)})^2}.
\end{align}
Similarly, we have
\begin{align*}
\frac{d}{dt}|_{t=t_2}\sqrt{b_{22}(t)}=\frac{|q_1|\dot{\alpha}(s_2)}{2\sqrt{\Delta(s_2)}}.
\end{align*}
So $\frac{d}{dt}|_{t=t_2}\sqrt{b_{22}(t)}=-1$ iff $q_1^2=4\Delta(s_2)\dot{\alpha}(s_2)^{-2}$. However,
\begin{eqnarray*}
\dot{\alpha}(s_2)&=&s_2^{-n}(-4(n+2)s_2^{n+1}+2p\kappa^{-1}s_2^n+c_1)\\
&=&2\kappa^{-1}(p-(2n+2)\kappa s_2-\frac{(p-(2n+2)\kappa s_1)(s_1^{n+1}-s_2^{n+1})}{(n+1)s_2^n(s_1-s_2)}),
\end{eqnarray*}
and
\begin{align*}
\Delta(s_2)=4\kappa^2s_2^{1-n}(\frac{p}{\kappa(2n+2)}-s_1)\frac{s_1^{n+1}-s_2^{n+1}}{s_1-s_2}.
\end{align*}
Thus
\begin{align}\label{q1}
q_1^2=\frac{2\kappa^3(p-(2n+2)\kappa s_1)(s_1^{n+1}-s_2^{n+1})}{(n+1)s_2^{n-1}(s_1-s_2)(p-(2n+2)\kappa s_2-\frac{(p-(2n+2)\kappa s_1)(s_1^{n+1}-s_2^{n+1})}{(n+1)s_2^n(s_1-s_2)})^2}.
\end{align}

We check now the identity $\psi^2=(q_1w_1+q_2w_2)^2$. Substituting \eqref{c1} and \eqref{c2f} into \eqref{ps} gives
\begin{eqnarray}
\label{p}\psi^2&=&4\kappa^{2n+2}\frac{s_1s_2}{(s_1-s_2)^2}(s_2^{n+1}-s_1^{n+1})(p-(2n+2)\kappa s_1)(p-(2n+2)\kappa s_2)\\
\notag&&\cdot(4(s_2^{n+1}-s_1^{n+1})-\frac{2p}{\kappa(n+1)}(s_2^n-s_1^n)).
\end{eqnarray}
On the other hand, by \eqref{w11} and \eqref{w22},
\begin{align*}
q_1w_1+q_2w_2=2q_1q_2\kappa^n((2n+2)(s_2^{n+1}-s_1^{n+1})-p\kappa^{-1}(s_2^n-s_1^n)).
\end{align*}
Thus
\begin{align}\label{qw}
(q_1w_1+q_2w_2)^2=4q_1^2q_2^2\kappa^{2n}((2n+2)(s_2^{n+1}-s_1^{n+1})-p\kappa^{-1}(s_2^n-s_1^n))^2.
\end{align}
We equate \eqref{p} with \eqref{qw}, and substitute $q_1^2$ and $q_2^2$ by \eqref{q2} and \eqref{q1}. The upshot is
\begin{eqnarray*}
\kappa^2&=&\frac{(n+1)s_1^ns_2^n}{2\kappa^2(s_2^{n+1}-s_1^{n+1})((2n+2)(s_2^{n+1}-s_1^{n+1})-p\kappa^{-1}(s_2^n-s_1^n))}\\
&&\cdot(p-(2n+2)\kappa s_1-\frac{(p-(2n+2)\kappa s_2)(s_2^{n+1}-s_1^{n+1})}{(n+1)s_1^n(s_2-s_1)})^2\\
&&\cdot(p-(2n+2)\kappa s_2-\frac{(p-(2n+2)\kappa s_1)(s_1^{n+1}-s_2^{n+1})}{(n+1)s_2^n(s_1-s_2)})^2.
\end{eqnarray*}
To simplify this expression, we introduce two quantities
\begin{eqnarray*}
x=(2n+2)p^{-1}\kappa s_1,&&y=(2n+2)p^{-1}\kappa s_2,
\end{eqnarray*}
i.e.,
\begin{eqnarray*}
s_1=\frac{px}{\kappa(2n+2)},&&s_2=\frac{py}{\kappa(2n+2)}.
\end{eqnarray*}
\begin{remark}
\eqref{s12} and \eqref{do} become $0<x<y<1$ and $y^{n+1}-y^n>x^{n+1}-x^n$.
\end{remark}
Now the formula for $\kappa^2$ takes the form
\begin{eqnarray*}
\kappa^2&=&\frac{(n+1)^2p^2x^ny^n}{(y^{n+1}-x^{n+1})(y^{n+1}-x^{n+1}-y^n+x^n)}\\
&&\cdot(1-x-\frac{(1-y)(y^{n+1}-x^{n+1})}{(n+1)x^n(y-x)})^2(1-y-\frac{(1-x)(x^{n+1}-y^{n+1})}{(n+1)y^n(x-y)})^2.
\end{eqnarray*}
Notice that $\kappa^2$ is determined by $x$ and $y$. Furthermore, \eqref{q2} and \eqref{q1} transform into
\begin{align}
\label{q11}&&q_1^2=\frac{p^2x^ny(1-x)}{(y-x)(y^{n+1}-x^{n+1}-y^n+x^n)}(1-x-\frac{(1-y)(y^{n+1}-x^{n+1})}{(n+1)x^n(y-x)})^2,\\
\label{q22}&&q_2^2=\frac{p^2xy^n(1-y)}{(y-x)(y^{n+1}-x^{n+1}-y^n+x^n)}(1-y-\frac{(1-x)(y^{n+1}-x^{n+1})}{(n+1)y^n(y-x)})^2.
\end{align}
It remains to show that
\begin{proposition}\label{cf}
Given a pair of non-zero integers $q_1$ and $q_2$ with $|q_1|>|q_2|$, there exists a point $(x,y)\in\Gamma=\{0<x<y<1;y^{n+1}-y^n>x^{n+1}-x^n\}$ such that \eqref{q11} and \eqref{q22} hold.
\end{proposition}
The proof of Proposition \ref{cf} needs several lemmas.
\begin{lemma}\label{yx}
For $(x,y)\in\Gamma$, $y^{k+1}-y^k>x^{k+1}-x^k$, $k=0,\cdots,n$.
\end{lemma}
\begin{proof}
Clearly the inequality is true for $k=n$. Assume the inequality holds when $k\ge m$. Then $y^{m+1}-x^{m+1}>y^m-x^m$, i.e., $\sum_{i=0}^my^{m-i}x^i>\sum_{j=0}^{m-1}y^{m-1-j}x^j$. So
\begin{eqnarray*}
\sum_{i=1}^my^{m-i}x^i&>&-y^m+y^{m-1}+\sum_{j=1}^{m-1}y^{m-1-j}x^j\\
&=&y^{m-1}(1-y)+\sum_{j=1}^{m-1}y^{m-1-j}x^j\\
&>&\sum_{j=1}^{m-1}y^{m-1-j}x^j,
\end{eqnarray*}
i.e., $\sum_{i=1}^my^{m-i}x^{i-1}>\sum_{j=1}^{m-1}y^{m-1-j}x^{j-1}$, or equivalently, $\sum_{i=0}^{m-1}y^{m-1-i}x^i>\sum_{j=0}^{m-2}y^{m-2-j}x^j$. Thus $y^m-x^m>y^{m-1}-x^{m-1}$, i.e., the inequality is true for $k=m-1$.
\end{proof}
Let
\begin{eqnarray*}
A_m=\sum_{i=0}^my^{m-i}x^i-\sum_{j=0}^{m-1}y^{m-1-j}x^j,&m=0,\cdots,n.
\end{eqnarray*}
For example, $A_0=1$, $A_1=y+x-1$. In the proof of Lemma \ref{yx}, we actually show that
\begin{corollary}\label{am}
$A_m>0$, $m=0,\cdots,n$.
\end{corollary}

\begin{lemma}\label{q12}
\begin{eqnarray*}
q_1^2=(\frac{p}{n+1})^2\frac{y(1-x)}{x^nA_n}(\sum_{i=0}^nx^iA_{n-i})^2,&&q_2^2=(\frac{p}{n+1})^2\frac{x(1-y)}{y^nA_n}(\sum_{i=0}^ny^iA_{n-i})^2.
\end{eqnarray*}
\end{lemma}
\begin{proof}
By \eqref{q11},
\begin{eqnarray*}
q_1^2&=&\frac{p^2x^ny(1-x)}{(y-x)^2A_n}(1-x-\frac{(1-y)\sum_{k=0}^nx^ky^{n-k}}{(n+1)x^n})^2\\
&=&\frac{p^2y(1-x)}{(n+1)^2x^n(y-x)^2A_n}(\sum_{k=0}^n(x^ky^{n-k+1}-x^{n+1}-x^ky^{n-k}+x^n))^2\\
&=&\frac{p^2y(1-x)}{(n+1)^2x^n(y-x)^2A_n}(\sum_{k=0}^nx^k(y^{n+1-k}-x^{n+1-k}-y^{n-k}+x^{n-k}))^2\\
&=&\frac{p^2y(1-x)}{(n+1)^2x^nA_n}(\sum_{k=0}^nx^kA_{n-k})^2
\end{eqnarray*}
By \eqref{q22}, we get the formula for $q_2^2$ in a similar way.
\end{proof}

\begin{lemma}
For $(x,y)\in\Gamma$, $q_1^2>q_2^2$.
\end{lemma}
\begin{proof}
By Lemma \ref{q12},
\begin{eqnarray*}
q_1^2-q_2^2&=&(\frac{p}{n+1})^2x^{-n}y^{-n}A_n^{-1}(y^{n+1}(1-x)(\sum_{i=0}^nx^iA_{n-i})^2-x^{n+1}(1-y)(\sum_{i=0}^ny^iA_{n-i})^2)\\
&=&(\frac{p}{n+1})^2x^{-n}y^{-n}A_n^{-1}((y^{n+1}-xy^{n+1})(\sum_{i=0}^nx^{2i}A_{n-i}^2+2\sum_{i<j}x^{i+j}A_{n-i}A_{n-j})\\
&&-(x^{n+1}-x^{n+1}y)(\sum_{i=0}^ny^{2i}A_{n-i}^2+2\sum_{i<j}y^{i+j}A_{n-i}A_{n-j}))\\
&=&(\frac{p}{n+1})^2x^{-n}y^{-n}A_n^{-1}(\sum_{i=0}^n(x^{2i}y^{n+1}-x^{2i+1}y^{n+1}-x^{n+1}y^{2i}+x^{n+1}y^{2i+1})A_{n-i}^2\\
&&+2\sum_{i<j}(x^{i+j}y^{n+1}-x^{i+j+1}y^{n+1}-x^{n+1}y^{i+j}+x^{n+1}y^{i+j+1})A_{n-i}A_{n-j}).\\
\end{eqnarray*}
\begin{claim}\label{pi}
$\Pi_k=x^ky^{n+1}-x^{k+1}y^{n+1}-x^{n+1}y^k+x^{n+1}y^{k+1}>0$, $k=0,\cdots,2n$.
\end{claim}
\begin{proof}
For $0\le k\le n$,
\begin{eqnarray*}
\Pi_k&=&x^ky^k(y^{n+1-k}-x^{n+1-k})-x^{k+1}y^{k+1}(y^{n-k}-x^{n-k})\\
&>&x^ky^k(y^{n+1-k}-x^{n+1-k}-y^{n-k}+x^{n-k})\\
&>&0.
\end{eqnarray*}
For $n+1\le k\le2n$,
\begin{eqnarray*}
\Pi_k=x^{n+1}y^{n+1}(y^{k-n}-x^{k-n}-y^{k-n-1}+x^{k-n-1})>0.
\end{eqnarray*}
\end{proof}
By Claim \ref{pi} and Corollary \ref{am}, $q_1^2-q_2^2>0$. This completes the proof.
\end{proof}

Let
\begin{eqnarray*}
G(z)=z^{n+1}-z^n.
\end{eqnarray*}
Then $G(0)=G(1)=0$, and
\begin{eqnarray*}
G'(z)=(n+1)z^{n-1}(z-\frac{n}{n+1})\left\{\begin{array}{ll}
<0,&z\in(0,\frac{n}{n+1})\\
>0,&z\in(\frac{n}{n+1},1).
\end{array}\right.
\end{eqnarray*}
Given $x\in[0,\frac{n}{n+1}]$, there exists a unique $y_\ast(x)\in[\frac{n}{n+1},1]$ such that $G(y_\ast(x))=G(x)$, i.e.,
\begin{align*}
y_\ast(x)^{n+1}-y_\ast(x)^n=x^{n+1}-x^n.
\end{align*}
Clearly, $y_\ast(x)$ is continuous and decreasing on $[0,\frac{n}{n+1}]$, with $y_\ast(0)=1$ and $y_\ast(\frac{n}{n+1})=\frac{n}{n+1}$.

Define
\begin{align*}
\Xi_1=\{(x,1):0<x\le1\},&&\Xi_2=\{(x,x):\frac{n}{n+1}<x\le1\},&&\Xi_3=\{(x,y_\ast(x)):0\le x\le\frac{n}{n+1}\}.
\end{align*}
Then $\partial\Gamma=\Xi_1\cup\Xi_2\cup\Xi_3$. By Lemma \ref{q12}, we have
\begin{lemma}\label{xi}
\begin{enumerate}
  \item On $\Xi_1$, $q_2^2\equiv0$.
  \item On $\Xi_2$, $q_1^2\equiv q_2^2$.
  \item On $\Xi_3$, $A_n\equiv0$.
\end{enumerate}
\end{lemma}
\begin{proof}[Proof of Proposition \ref{cf}]
For $0\le a<1$, let us consider the line $c_a(\tau)=(\tau,\tau+a)$, $\tau\in\mathbb{R}$. By Lemma \ref{xi}.1, $q_2^2=0$ at $c_a(1-a)=(1-a,1)$; when $\tau$ decreases, $c_a(\tau)$ approaches $\Xi_3$, and hence $q_2^2\rightarrow\infty$ as $A_n\rightarrow0$ (cf. Lemma \ref{xi}.3). Thus given $L_2\in\mathbb{N}$, there exists the largest real number $\Theta(L_2;a)<1-a$ such that $q_2^2=L_2^2$ at $c_a(\Theta(L_2;a))=(\Theta(L_2;a),\Theta(L_2;a)+a)$. Notice that $\Theta(L_2;a)$ is a continuous function of $a$, and $\Theta(L_2;a)\rightarrow0^+$ when $a\rightarrow1^-$.

We study now the behavior of $q_1^2$ along the segment $\Theta(L_2;a)$, $0\le a<1$. By Lemma \ref{q12},
\begin{eqnarray*}
\frac{q_1^2}{q_2^2}&=&\frac{y^{n+1}(1-x)(\sum_{i=0}^nx^iA_{n-i})^2}{x^{n+1}(1-y)(\sum_{i=0}^ny^iA_{n-i})^2}\\
&>&\frac{y^{n+1}(1-x)\cdot x^{2n}}{x^{n+1}(1-y)(\sum_{i=0}^ny^iA_{n-i})^2}\\
&=&\frac{y^{n+1}(1-x)}{(\sum_{i=0}^ny^iA_{n-i})^2}\cdot\frac{x^{n-1}}{1-y}.
\end{eqnarray*}
Since $y^n-y^{n+1}<x^n-x^{n+1}$, we have $y^n(1-y)<x^n(1-x)$, i.e., $1-y<y^{-n}x^n(1-x)$. Thus
\begin{eqnarray*}
\frac{q_1^2}{q_2^2}&>&\frac{y^{n+1}(1-x)}{(\sum_{i=0}^ny^iA_{n-i})^2}\cdot\frac{x^{n-1}y^n}{x^n(1-x)}=\frac{y^{2n+1}}{(\sum_{i=0}^ny^iA_{n-i})^2}\cdot\frac{1}{x}\rightarrow\infty,\;\;\mbox{as }x\rightarrow0\mbox{ and }y\rightarrow1.
\end{eqnarray*}
When $a=0$, $c_0(\Theta(L_2;0))=(\Theta(L_2;0),\Theta(L_2;0))$ lies on $\Xi_2$. By Lemma \ref{xi}.2, $q_1^2=q_2^2=L_2^2$ at $c_0(\Theta(L_2;0))$. On the other hand, $\Theta(L_2;a)$ is getting closer to $0$ when $a$ increases to $1$. Thus given $L_1\in\mathbb{N}$ with $L_1>L_2$, the above estimate asserts an $a(L_1;L_2)\in(0,1)$ such that $q_1^2=L_1^2$ at $c_{a(L_1;L_2)}(\Theta(L_2;a(L_1;L_2)))$.

In summary, given a pair of integers $L_1>L_2>0$, let $q_i^2$ be as in Lemma \ref{q12} with $x=\Theta(L_2;a(L_1;L_2))$ and $y=\Theta(L_2;a(L_1;L_2))+a(L_1;L_2)$. Then $q_i^2=L_i^2$. This completes the proof.
\end{proof}

\begin{remark}
As before, it is clear that $b_{ij}(t)$ and $c(t)$ are smooth and even in $t$ around $t_k$, $k=1,2$.
\end{remark}
By now, we have checked that all smooth collapse conditions are satisfied, and hence construct a complete positive Einstein metric on the $3$-sphere bundle associated with $P_q$. This finishes the proof of Theorem \ref{pos}.

\begin{remark}\label{vv}
The Einstein manifold constructed above has volume $C\kappa^n(n+1)^{-1}(s_2^{n+1}-s_1^{n+1})$ (cf. \S\ref{aav}). As pointed out in Remark \ref{v}, it is hopeless to use this formula to verify whether two Einstein manifolds have the same volume.
\end{remark}

\subsection{Homeomorphism and diffeomorphism classifications}

We specialize now to the case that the base $V$ is the complex projective plane $\mathbb{C}P^2$. Let $P_q$ be the principal $2$-torus bundle over $\mathbb{C}P^2$ with characteristic classes $(q_1a,q_2a)$, where $a$ is a generator of $H^2(\mathbb{C}P^2;\mathbb{Z})\cong\mathbb{Z}$. Let $E_1=P_q\times_{S^1\times S^1}S^3$, $E_2=P_q\times_{S^1\times S^1}D^4$, and $E_3=P_q\times_{S^1\times S^1}\mathbb{C}^2$ be respectively the standard $3$-sphere bundle, $4$-ball bundle, and complex plane bundle associated with $P_q$. Denote by $\pi_i:E_i\rightarrow\mathbb{C}P^2$ the projection of bundle. There are natural inclusions $j_{kl}:E_k\hookrightarrow E_{l}$, $k<l$; we have $j_{13}=j_{23}\circ j_{12}$.

\subsubsection{Characteristic classes}

We calculate here the second Stiefel-Whitney class and the first Pontrjagin class of tangent bundle $TE_i$. But the starting point is the Chern classes of $E_3$.
\begin{lemma}[Chern class]\label{cc}
$c_1(E_3)=(q_1+q_2)a$ and $c_2(E_3)=q_1q_2a^2$.
\end{lemma}
\begin{proof}
There is a natural splitting $E_3\cong L_{q_1}\oplus L_{q_2}$, where $L_{q_i}$ is the complex line bundle over $\mathbb{C}P^2$ with first Chern class $q_ia$. So the total Chern class of $E_3$ is
\begin{eqnarray*}
c(E_3)=c(L_{q_1})\cup c(L_{q_2})=(1+q_1a)\cup(1+q_2a)=1+(q_1+q_2)a+q_1q_2a^2.
\end{eqnarray*}
\end{proof}

It follows from the preceding lemma that the $3$-sphere bundle $E_1$ has Euler class $q_1q_2a^2$. Applying to $E_1$ the Gysin sequence shows that
\begin{lemma}\label{h}
\begin{align*}
\left.\begin{array}{llll}
H^0(E_1;\mathbb{Z})=\mathbb{Z}&H^1(E_1;\mathbb{Z})=0&H^2(E_1;\mathbb{Z})=\mathbb{Z}\cdot\pi_1^\ast a&H^3(E_1;\mathbb{Z})=0\\
H^4(E_1;\mathbb{Z})=\mathbb{Z}_{|q_1q_2|}\cdot(\pi_1^\ast a)^2&H^5(E_1;\mathbb{Z})=\mathbb{Z}\cdot x&H^6(E_1;\mathbb{Z})=0&H^7(E_1;\mathbb{Z})=\mathbb{Z}\cdot \pi_1^\ast a\cup x
\end{array}\right.
\end{align*}
\end{lemma}

\begin{lemma}[Second Stiefel-Whitney class]\label{w}
$w_2(TE_i)=(1+q_1+q_2)\pi_i^\ast a$ (mod $2$).
\end{lemma}
\begin{proof}
Since $TE_3\cong\pi_3^\ast T\mathbb{C}P^2\oplus\pi_3^\ast E_3$, the total Stiefel-Whitney class of $TE_3$ is
\begin{eqnarray*}
w(TE_3)=\pi_3^\ast(w(T\mathbb{C}P^2)\cup w(E_3)).
\end{eqnarray*}
By Lemma \ref{cc},
\begin{eqnarray*}
w_2(TE_3)=\pi_3^\ast(w_2(T\mathbb{C}P^2)+w_2(E_3))=\pi_3^\ast(3a+(q_1+q_2)a)=(1+q_1+q_2)\pi_3^\ast a&&(\mbox{mod }2).
\end{eqnarray*}
Since $TE_2\cong j_{23}^\ast TE_3$, we have
\begin{eqnarray*}
w_2(TE_2)=j_{23}^\ast w_2(TE_3)=(1+q_1+q_2)(\pi_3\circ j_{23})^\ast a=(1+q_1+q_2)\pi_2^\ast a&&(\mbox{mod }2).
\end{eqnarray*}
Finally, to determine $w_2(TE_1)$, we need first to compute $w_2(\ker(d\pi_1))$ as $TE_1\cong\pi_1^\ast T\mathbb{C}P^2\oplus\ker(d\pi_1)$. Notice that $\pi_1^\ast E_3\cong\ker(d\pi_1)\oplus\nu^1$, where $\nu^1$ is the rank $1$ normal bundle of $E_1$ in $E_3$. So $w_2(\pi_1^\ast E_3)=w_2(\ker(d\pi_1))$, i.e., $w_2(\ker(d\pi_1))=\pi_1^\ast w_2(E_3)=(q_1+q_2)\pi_1^\ast a$ (mod $2$) (cf. Lemma \ref{cc}). Thus
\begin{eqnarray*}
w_2(TE_1)=\pi_1^\ast w_2(T\mathbb{C}P^2)+w_2(\ker(d\pi_1))=\pi_1^\ast(3a)+(q_1+q_2)\pi_1^\ast a=(1+q_1+q_2)\pi_1^\ast a&&(\mbox{mod }2).
\end{eqnarray*}
\end{proof}

\begin{corollary}
$E_i$ is spin iff $q_1+q_2$ is odd.
\end{corollary}

\begin{lemma}[First Pontrjagin class]\label{fp}
$p_1(TE_i)=(3+q_1^2+q_2^2)(\pi_i^\ast a)^2$, $i=2,3$.
\end{lemma}
\begin{proof}
By \cite[Theorem 15.3]{MilSta74}, $2(p(TE_3)-p(\pi_3^\ast T\mathbb{C}P^2)\cup p(\pi_3^\ast E_3))=0$. Since $H^\ast(E_3;\mathbb{Z})\cong H^\ast(\mathbb{C}P^2;\mathbb{Z})$ is torsion-free, we have $p(TE_3)=\pi_3^\ast(p(T\mathbb{C}P^2)\cup p(E_3))$. By \cite[Corollary 15.5]{MilSta74} and Lemma \ref{cc}, $p_1(E_3)=c_1(E_3)^2-2c_2(E_3)=(q_1^2+q_2^2)a^2$. Therefore
\begin{eqnarray*}
p_1(TE_3)=\pi_3^\ast(p_1(T\mathbb{C}P^2)+p_1(E_3))=\pi_3^\ast(3a^2+(q_1^2+q_2^2)a^2)=(3+q_1^2+q_2^2)(\pi_3^\ast a)^2.
\end{eqnarray*}
Moreover, $p_1(TE_2)=j_{23}^\ast p_1(TE_3)=(3+q_1^2+q_2^2)(\pi_2^\ast a)^2$.
\end{proof}

In addition to $j_{23}:E_2\rightarrow E_3$, there are three more inclusions. The first one is $j_{30}:E_3\rightarrow(E_3,E_3^0)$, where $E_3^0$ is the set of all non-zero vectors in $E_3$. The second one is $j_{21}:E_2\rightarrow(E_2,E_1)$, which induces a monomorphism $j_{21}^\ast:H^4(E_2,E_1;\mathbb{Z})\rightarrow H^4(E_2;\mathbb{Z})$ since $H^3(E_1;\mathbb{Z})=0$. The third one is $j_{10}:(E_2,E_1)\rightarrow(E_3,E_3^0)$, which induces an isomorphism $j_{10}^\ast:H^4(E_3,E_3^0;\mathbb{Z})\rightarrow H^4(E_2,E_1;\mathbb{Z})$.
\begin{remark}\label{j}
We have $j_{10}\circ j_{21}=j_{30}\circ j_{23}:E_2\rightarrow(E_3,E_3^0)$.
\end{remark}
\begin{lemma}\label{t}
Let $U\in H^4(E_3,E_3^0;\mathbb{Z})$ be the Thom class. Then $j_{21}^\ast j_{10}^\ast U=q_1q_2(\pi_2^\ast a)^2$.
\end{lemma}
\begin{proof}
By definition, $E_3$ has Euler class $e(E_3)=(\pi_3^\ast)^{-1}(j_{30}^\ast U)$. So $(j_{30}\circ j_{23})^\ast U=(\pi_3\circ j_{23})^\ast(e(E_3))$. By Remark \ref{j}, we have $(j_{10}\circ j_{21})^\ast U=\pi_2^\ast(e(E_3))$. The lemma follows now from Lemma \ref{cc}.
\end{proof}

\subsubsection{Kreck-Stolz invariants} The $3$-sphere bundle $E_1$ is a simply-connected closed $7$-manifold with integral cohomology ring of the type $H^2\cong\mathbb{Z}$, $H^3=0$, and $H^4$ being a cyclic group of finite order generated by the square of a generator of $H^2$ (cf. Lemma \ref{h}). Kreck and Stolz \cite{KreSto88, KreSto91} were able to define a triple of $\mathbb{Q}/\mathbb{Z}$-valued invariants to determine the homeomorphism and diffeomorphism types of this class of manifolds as follows.

Given a smooth $7$-manifold $X$ as above, let $x$ be a generator of $H^2(X;\mathbb{Z})$. By \cite[Proposition 2.2]{KreSto91}, we can find an $8$-manifold $Y$ with $y\in H^2(Y;\mathbb{Z})$ such that $\partial Y=X$, $y$ restricts to $x$ on $X$, $w_2(Y)=0$ when $X$ is spin, and $w_2(Y)=y$ (mod $2$) when $X$ is non-spin. Since $H^3(X;\mathbb{Q})=H^4(X;\mathbb{Q})=0$, we have $H^4(Y,X;\mathbb{Q})\cong H^4(Y;\mathbb{Q})$. Via this isomorphism, we can regard $y^2$ and $p_1(TY)$ as elements of $H^4(Y,X;\mathbb{Q})$, and evaluate their products on the relative fundamental class $[Y,X]$. By abuse of notation,, we write $y^4=\langle y^4,[Y,X]\rangle$, $p_1^2=\langle p_1(TY)^2,[Y,X]\rangle$, and $y^2p_1=\langle y^2\cup p_1(TY),[Y,X]\rangle$. Denote by $\mbox{sign}(Y)$ the signature of $Y$.
\begin{definition}
\begin{enumerate}
  \item In the spin case, define
\begin{align*}
s_1(Y,y)=\frac{p_1^2}{896}-\frac{\mbox{sign}(Y)}{224},&&s_2(Y,y)=\frac{y^4}{24}-\frac{y^2p_1}{48},&&s_3(Y,y)=\frac{2}{3}y^4-\frac{y^2p_1}{12}.
\end{align*}
  \item In the non-spin case, define
\begin{align*}
s_1(Y,y)=\frac{y^4}{384}-\frac{y^2p_1}{192}+\frac{p_1^2}{896}-\frac{\mbox{sign}(Y)}{224},&&s_2(Y,y)=\frac{5}{24}y^4-\frac{y^2p_1}{24},&&s_3(Y,y)=\frac{13}{8}y^4-\frac{y^2p_1}{8}.
\end{align*}
\end{enumerate}
\end{definition}
By \cite[Proposition 2.5]{KreSto91}, $s_i(Y,y)$ (mod $\mathbb{Z}$) depends only on $(X,x)$. So we define $s_i(X,x)=s_i(Y,y)$ (mod $\mathbb{Z}$). Moreover, $s_i(X,x)$ is independent of the choice of a generator of $H^2(X;\mathbb{Z})$. This leads to
\begin{definition}[Kreck-Stolz invariants]
$s_i(X)=s_i(X,x)$, $i=1,2,3$.
\end{definition}
The next theorem is due to Kreck and Stolz \cite{KreSto88, KreSto91} (see \cite[Theorem 1]{KreSto98} for the correct formulation).
\begin{theorem}\label{kss}
Let $X_1$ and $X_2$ be smooth $7$-manifolds of the above type such that $|H^4(X_1;\mathbb{Z})|=|H^4(X_2;\mathbb{Z})|$ which are both spin or both non-spin. Then $X_1$ is diffeomorphic (homeomorphic) to $X_2$ iff $s_i(X_1)=s_i(X_2)$ for $i=1,2,3$ (resp. $28s_1(X_1)=28s_1(X_2)$ and $s_i(X_1)=s_i(X_2)$ for $i=2,3$).
\end{theorem}

\subsubsection{Calculation of Kreck-Stolz invariants associated with $E_1$} By Lemma \ref{h}, $x=\pi_1^\ast a$ is a generator of $H^2(E_1;\mathbb{Z})$.
\begin{lemma}
$E_2$, together with $y=\pi_2^\ast a$, qualifies as the required zero bordism.
\end{lemma}
\begin{proof}
Clearly $x=j_{12}^\ast y$ for the inclusion $j_{12}:E_1\rightarrow E_2$. The lemma follows now from Lemma \ref{w}.
\end{proof}

\begin{remark}
With rational coefficients, $j_{21}^\ast:H^4(E_2,E_1;\mathbb{Q})\rightarrow H^4(E_2;\mathbb{Q})\cong\mathbb{Q}$ is an isomorphism. By Lemma \ref{t}, the image of the Thom class in $H^4(E_2,E_1;\mathbb{Q})$ is $j_{10}^\ast U=q_1q_2(j_{21}^\ast)^{-1}(y^2)$.
\end{remark}

\begin{lemma}
The involved relative characteristic numbers are
\begin{align*}
y^4=\frac{1}{q_1q_2},&&y^2p_1=\frac{3+q_1^2+q_2^2}{q_1q_2},&&p_1^2=\frac{(3+q_1^2+q_2^2)^2}{q_1q_2}.
\end{align*}
\end{lemma}
\begin{proof}
Since $(j_{21}^\ast)^{-1}(y^2)=(q_1q_2)^{-1}j_{10}^\ast U$, we have
\begin{eqnarray*}
y^4&=&\langle(j_{21}^\ast)^{-1}(y^2)\cup(j_{21}^\ast)^{-1}(y^2),[E_2,E_1]\rangle\\
&=&\frac{1}{q_1q_2}\langle(j_{21}^\ast)^{-1}(y^2)\cup j_{10}^\ast U,[E_2,E_1]\rangle\\
&=&\frac{1}{q_1q_2}\langle a^2,[\mathbb{C}P^2]\rangle\\
&=&\frac{1}{q_1q_2}.
\end{eqnarray*}
By Lemma \ref{fp}, $p_1(TE_2)=(3+q_1^2+q_2^2)y^2$. So
\begin{align*}
p_1^2=(3+q_1^2+q_2^2)^2y^4=\frac{(3+q_1^2+q_2^2)^2}{q_1q_2},&&y^2p_1=(3+q_1^2+q_2^2)y^4=\frac{3+q_1^2+q_2^2}{q_1q_2}.
\end{align*}
\end{proof}
\begin{corollary}
$\mbox{sign}(E_2)=\mbox{sgn}(q_1q_2)$.
\end{corollary}
\begin{corollary}\label{kse}
The Kreck-Stolz invariants associated with $E_1$ are
\begin{enumerate}
  \item In the spin case, i.e., when $q_1+q_2$ is odd,
\begin{align*}
s_1(E_1)=\frac{(3+q_1^2+q_2^2)^2}{896q_1q_2}-\frac{\mbox{sgn}(q_1q_2)}{224},
&&s_2(E_1)=-\frac{1+q_1^2+q_2^2}{48q_1q_2},
&&s_3(E_1)=\frac{5-q_1^2-q_2^2}{12q_1q_2}.
\end{align*}
  \item In the non-spin case, i.e., when $q_1+q_2$ is even,
\begin{align*}
s_1(E_1)=\frac{1}{384q_1q_2}-\frac{3+q_1^2+q_2^2}{192q_1q_2}+\frac{(3+q_1^2+q_2^2)^2}{896q_1q_2}-\frac{\mbox{sgn}(q_1q_2)}{224},&&s_2(E_1)=\frac{2-q_1^2-q_2^2}{24q_1q_2},\\
s_3(E_1)=\frac{10-q_1^2-q_2^2}{8q_1q_2}.
\end{align*}
\end{enumerate}
\end{corollary}

\subsubsection{Classifications} Assume that $E_1$ and $\widehat{E}_1$ are respectively the $3$-sphere bundles associated with $P_q$ and $P_{\widehat{q}}$. By Theorem \ref{kss} and Lemma \ref{h}, a necessary condition for $E_1$ homeomorphic to $\widehat{E}_1$ is $|q_1q_2|=|\widehat{q}_1\widehat{q}_2|$.
\begin{proposition}\label{m}
Assume $q_1q_2=-\widehat{q}_1\widehat{q}_2$. Then $E_1$ is homeomorphic (diffeomorphic) to $\widehat{E}_1$ iff $|q_1q_2|=1$.
\end{proposition}
\begin{proof}
Suppose $q_1q_2=-\widehat{q}_1\widehat{q}_2=1$. So $q_1^2+q_2^2=\widehat{q}_1^2+\widehat{q}_2^2=2$. By Corollary \ref{kse}, $s_i(E_1)-s_i(\widehat{E}_1)\in\mathbb{Z}$, $i=1,2,3$. By Theorem \ref{kss}, $E_1$ is diffeomorphic, and hence homeomorphic, to $\widehat{E}_1$. Conversely, suppose $M=q_1q_2=-\widehat{q}_1\widehat{q}_2$ and $E_1$ is homeomorphic to $\widehat{E}_1$. There are two cases.

When $E_1$ and $\widehat{E}_1$ are spin, i.e., $q_1+q_2$ and $\widehat{q}_1+\widehat{q}_2$ are odd, we see that $M$ is even. By Theorem \ref{kss} and Corollary \ref{kse}.1,
\begin{align*}
s_2(E_1)-s_2(\widehat{E}_1)=-\frac{2+q_1^2+q_2^2+\widehat{q}_1^2+\widehat{q}_2^2}{48M}\in\mathbb{Z},&&s_3(E_1)-s_3(\widehat{E}_1)=\frac{10-q_1^2-q_2^2-\widehat{q}_1^2-\widehat{q}_2^2}{12M}\in\mathbb{Z},
\end{align*}
i.e., $2+q_1^2+q_2^2+\widehat{q}_1^2+\widehat{q}_2^2=48Mk$ and $10-q_1^2-q_2^2-\widehat{q}_1^2-\widehat{q}_2^2=12Ml$ for $k,l\in\mathbb{Z}$. So $12=48Mk+12Ml$, i.e., $1=M(4k+l)$, and hence $M=\pm1$. This gives a contradiction as $M$ is even.

When $E_1$ and $\widehat{E}_1$ are non-spin, i.e., $q_1+q_2$ and $\widehat{q}_1+\widehat{q}_2$ are even, it follows from Theorem \ref{kss} and Corollary \ref{kse}.2 that
\begin{align*}
s_2(E_1)-s_2(\widehat{E}_1)=\frac{4-q_1^2-q_2^2-\widehat{q}_1^2-\widehat{q}_2^2}{24M}\in\mathbb{Z},&&s_3(E_1)-s_3(\widehat{E}_1)=\frac{20-q_1^2-q_2^2-\widehat{q}_1^2-\widehat{q}_2^2}{8M}\in\mathbb{Z},
\end{align*}
i.e., $4-q_1^2-q_2^2-\widehat{q}_1^2-\widehat{q}_2^2=24M\widehat{k}$ and $20-q_1^2-q_2^2-\widehat{q}_1^2-\widehat{q}_2^2=8M\widehat{l}$ for $\widehat{k},\widehat{l}\in\mathbb{Z}$.
So $16=8M\widehat{l}-24M\widehat{k}$, i.e., $2=M(\widehat{l}-3\widehat{k})$, and hence $M=\pm1$ or $\pm2$. If $M=\pm2$, then $q_1$ and $q_2$ have to be even. So $M$ is divisible by $4$. This gives a contradiction. If $M=\pm1$, then $q_1^2+q_2^2=\widehat{q}_1^2+\widehat{q}_2^2=2$. By Corollary \ref{kse}.2, $s_1(E_1)=s_1(\widehat{E}_1)=0$. By Theorem \ref{kss}, $E_1$ is diffeomorphic to $\widehat{E}_1$.
\end{proof}

From now on, we assume $M=q_1q_2=\widehat{q}_1\widehat{q}_2$. Let $N=q_1^2+q_2^2$ and $\widehat{N}=\widehat{q}_1^2+\widehat{q}_2^2$.
\begin{lemma}
$q_1+q_2\equiv q_1^2+q_2^2$ (mod $2$).
\end{lemma}
\begin{proof}
Clearly $q_1+q_2\equiv(q_1+q_2)^2=q_1^2+q_2^2+2q_1q_2\equiv q_1^2+q_2^2$ (mod $2$).
\end{proof}

Let us consider first the spin case, i.e., $N\equiv\widehat{N}\equiv1$ (mod $2$). By Corollary \ref{kse}.1,
\begin{align*}
s_1(E_1)-s_1(\widehat{E}_1)=\frac{(6+N+\widehat{N})(N-\widehat{N})}{896M},&&s_2(E_1)-s_2(\widehat{E}_1)=\frac{\widehat{N}-N}{48M},&&s_3(E_1)-s_3(\widehat{E}_1)=\frac{\widehat{N}-N}{12M}.
\end{align*}

\begin{proposition}\label{sp}
In the spin case, $E_1$ is homeomorphic (diffeomorphic) to $\widehat{E}_1$ iff $N\equiv\widehat{N}$ (mod $48M$) (and $N+(\frac{N+1}{2})^2\equiv\widehat{N}+(\frac{\widehat{N}+1}{2})^2$ (mod $224M$)).
\end{proposition}
\begin{proof}
If $E_1$ is homeomorphic to $\widehat{E}_1$, then $s_2(E_1)-s_2(\widehat{E}_1)\in\mathbb{Z}$, i.e., $\widehat{N}\equiv N$ (mod $48M$). Conversely, if $\widehat{N}\equiv N$ (mod $48M$), i.e., $\widehat{N}-N=48Mk$ for some $k\in\mathbb{Z}$, then $s_2(E_1)-s_2(\widehat{E}_1)\in\mathbb{Z}$ and $s_3(E_1)-s_3(\widehat{E}_1)\in\mathbb{Z}$. Moreover, $28(s_1(E_1)-s_1(\widehat{E}_1))=-3k(3+\frac{N+\widehat{N}}{2})\in\mathbb{Z}$ as $N\equiv\widehat{N}$ (mod $2$). By Theorem \ref{kss}, $E_1$ is homeomorphic to $\widehat{E}_1$.

For $E_1$ diffeomorphic to $\widehat{E}_1$, we need furthermore
\begin{align*}
s_1(E_1)-s_1(\widehat{E}_1)=\frac{N+(\frac{N+1}{2})^2-\widehat{N}-(\frac{\widehat{N}+1}{2})^2}{224M}\in\mathbb{Z},
\end{align*}
which is equivalent to $N+(\frac{N+1}{2})^2\equiv\widehat{N}+(\frac{\widehat{N}+1}{2})^2$ (mod $224M$).
\end{proof}

We turn now to the non-spin case, i.e., $N\equiv\widehat{N}\equiv0$ (mod $2$). By Corollary \ref{kse}.2,
\begin{align*}
s_1(E_1)-s_1(\widehat{E}_1)=\frac{(6+N+\widehat{N})(N-\widehat{N})}{896M}+\frac{\widehat{N}-N}{192M},&&s_2(E_1)-s_2(\widehat{E}_1)=\frac{\widehat{N}-N}{24M},\\
s_3(E_1)-s_3(\widehat{E}_1)=\frac{\widehat{N}-N}{8M}.
\end{align*}

\begin{proposition}\label{nsp}
In the non-spin case, $E_1$ is homeomorphic (diffeomorphic) to $\widehat{E}_1$ iff $N\equiv\widehat{N}$ (mod $24M$) (and $N+3(\frac{N}{2})^2\equiv\widehat{N}+3(\frac{\widehat{N}}{2})^2$ (mod $672M$)).
\end{proposition}
\begin{proof}
If $E_1$ is homeomorphic to $\widehat{E}_1$, then $s_2(E_1)-s_2(\widehat{E}_1)\in\mathbb{Z}$, i.e., $\widehat{N}\equiv N$ (mod $24M$). Conversely, if $\widehat{N}\equiv N$ (mod $24M$), i.e., $\widehat{N}-N=24Mk$ for some $k\in\mathbb{Z}$, then $s_2(E_1)-s_2(\widehat{E}_1)\in\mathbb{Z}$ and $s_3(E_1)-s_3(\widehat{E}_1)\in\mathbb{Z}$. Moreover, $28(s_1(E_1)-s_1(\widehat{E}_1))=-k(1+\frac{3}{2}N+18Mk)\in\mathbb{Z}$ as $N$ is even. By Theorem \ref{kss}, $E_1$ is homeomorphic to $\widehat{E}_1$.

For $E_1$ diffeomorphic to $\widehat{E}_1$, we need furthermore
\begin{align*}
s_1(E_1)-s_1(\widehat{E}_1)=\frac{N+3(\frac{N}{2})^2-\widehat{N}-3(\frac{\widehat{N}}{2})^2}{672M}\in\mathbb{Z},
\end{align*}
which is equivalent to $N+3(\frac{N}{2})^2\equiv\widehat{N}+3(\frac{\widehat{N}}{2})^2$ (mod $672M$).
\end{proof}

\begin{proof}[Proof of Theorem \ref{cla}]
Theorem \ref{cla} follows from Propositions \ref{m}, \ref{sp}, and \ref{nsp}.
\end{proof}

\section*{Acknowledgements}
This work forms part of the author's PhD thesis. He is grateful to his supervisor, Prof. M. Y. Wang, for his consistent help and inspiration underlying all of the work. He wishes to thank Drs. J. Cuadros, A. S. Dancer, C. R. Graham, P. Guan, M. Min-Oo, A. J. Nicas and the anonymous referee for their interests in this paper and many very valuable comments.

\bibliographystyle{amsplain}
\bibliography{DCbib}

\providecommand{\bysame}{\leavevmode\hbox to3em{\hrulefill}\thinspace}
\providecommand{\MR}{\relax\ifhmode\unskip\space\fi MR }
\providecommand{\MRhref}[2]{%
  \href{http://www.ams.org/mathscinet-getitem?mr=#1}{#2}
}
\providecommand{\href}[2]{#2}
\begin{thebibliography}{10}

\bibitem{And05}
M.~T. Anderson, \emph{Geometric aspects of the {A}d{S}/{CFT} correspondence},
  AdS/CFT correspondence: Einstein metrics and their conformal boundaries, IRMA
  Lect. Math. Theor. Phys., vol.~8, Eur. Math. Soc., Z\"{u}rich, 2005,
  pp.~1--31.

\bibitem{And08}
\bysame, \emph{Einstein metrics with prescribed conformal infinity on
  4-manifolds}, Geom. Funct. Anal. \textbf{18} (2008), no.~2, 305--366.

\bibitem{Ber82}
L.~B\'{e}rard-Bergery, \emph{Sur de nouvelles vari\'{e}t\'{e}s riemanniennes
  d'{E}instein}, Inst. \'{E}lie Cartan, vol.~6, Univ. Nancy, Nancy, 1982,
  pp.~1--60.

\bibitem{Bes87}
A.~L. Besse, \emph{Einstein manifolds}, Ergebnisse der Mathematik und ihrer
  Grenzgebiete (3), vol.~10, Springer-Verlag, Berlin, 1987.

\bibitem{BohWanZil04}
C.~B\"{o}hm, M.~Wang, and W.~Ziller, \emph{A variational approach for compact
  homogeneous {E}instein manifolds}, Geom. Funct. Anal. \textbf{14} (2004),
  no.~4, 681--733.

\bibitem{Bre03}
S.~Brendle, \emph{Global existence and convergence for a higher order flow in
  conformal geometry}, Ann. of Math. (2) \textbf{158} (2003), no.~1, 323--343.

\bibitem{ChaYan95}
S.-Y.~A. Chang and P.~C. Yang, \emph{Extremal metrics of zeta function
  determinants on $4$-manifolds}, Ann. of Math. (2) \textbf{142} (1995), no.~1,
  171--212.

\bibitem{Che09}
D.~Chen, \emph{Construction of conformally compact {E}instein manifolds},
  arXiv:0908.1430v2, 2009.

\bibitem{DjaMal08}
Z.~Djadli and A.~Malchiodi, \emph{Existence of conformal metrics with constant
  ${Q}$-curvature}, Ann. of Math. (2) \textbf{168} (2008), no.~3, 813--858.

\bibitem{EscWan00}
J.-H. Eschenburg and M.~Y. Wang, \emph{The initial value problem for
  cohomogeneity one {E}instein metrics}, J. Geom. Anal. \textbf{10} (2000),
  no.~1, 109--137.

\bibitem{FefGra02}
C.~Fefferman and C.~R. Graham, \emph{Q-curvature and {P}oincar\'{e} metrics},
  Math. Res. Lett. \textbf{9} (2002), 139--151.

\bibitem{GibHarYas04}
G.~W. Gibbons, S.~A. Hartnoll, and Y.~Yasui, \emph{Properties of some
  five-dimensional {E}instein metrics}, Classical Quantum Gravity \textbf{21}
  (2004), no.~19, 4697--4730.

\bibitem{Gra00}
C.~R. Graham, \emph{Volume and area renormalizations for conformally compact
  {E}instein metrics}, The Proceedings of the 19th Winter School ``Geometry and
  Physics" (Srn\'{\i}, 1999), Rend. Circ. Mat. Palermo (2) Suppl., no.~63,
  2000, pp.~31--42.

\bibitem{GraZwo03}
C.~R. Graham and M.~Zworski, \emph{Scattering matrix in conformal geometry},
  Invent. Math. \textbf{152} (2003), no.~1, 89--118.

\bibitem{HasSakYas05}
Y.~Hashimoto, M.~Sakaguchi, and Y.~Yasui, \emph{New infinite series of
  {E}instein metrics on sphere bundles from {A}d{S} black holes}, Comm. Math.
  Phys. \textbf{257} (2005), no.~2, 273--285.

\bibitem{KoiSak86}
N.~Koiso and Y.~Sakane, \emph{Nonhomogeneous {K}\"{a}hler-{E}instein metrics on
  compact complex manifolds}, Curvature and topology of Riemannian manifolds
  (Katata, 1985), Lecture Notes in Math., no. 1201, Springer, Berlin, 1986,
  pp.~165--179.

\bibitem{KreSto88}
M.~Kreck and S.~Stolz, \emph{A diffeomorphism classification of $7$-dimensional
  homogeneous {E}instein manifolds with
  ${SU}(3)\times{SU}(2)\times{U}(1)$-symmetry}, Ann. of Math. (2) \textbf{127}
  (1988), no.~2, 373--388.

\bibitem{KreSto91}
\bysame, \emph{Some nondiffeomorphic homeomorphic homogeneous $7$-manifolds
  with positive sectional curvature}, J. Differential Geom. \textbf{33} (1991),
  no.~2, 465--486.

\bibitem{KreSto98}
\bysame, \emph{A correction on: ``{S}ome nondiffeomorphic homeomorphic
  homogeneous $7$-manifolds with positive sectional curvature"}, J.
  Differential Geom. \textbf{49} (1998), no.~1, 203--204.

\bibitem{Kru05}
B.~Kruggel, \emph{Homeomorphism and diffeomorphism classification of
  {E}schenburg spaces}, Q. J. Math. \textbf{56} (2005), no.~4, 553--577.

\bibitem{Lau02}
J.~Lauret, \emph{Finding {E}instein solvmanifolds by a variational method},
  Math. Z. \textbf{241} (2002), no.~1, 83--99.

\bibitem{LotShe00}
J.~Lott and Z.~Shen, \emph{Manifolds with quadratic curvature decay and slow
  volume growth}, Ann. Sci. \'{E}cole Norm. Sup. (4) \textbf{33} (2000), no.~2,
  275--290.

\bibitem{LvPagPop04}
H.~L\"{u}, D.~N. Page, and C.~N. Pope, \emph{New inhomogeneous {E}instein
  metrics on sphere bundles over {E}instein-{K}\"{a}hler manifolds}, Phys.
  Lett. B \textbf{593} (2004), no.~1-4, 218--226.

\bibitem{MilSta74}
J.~W. Milnor and J.~D. Stasheff, \emph{Characteristic classes}, Annals of
  Mathematics Studies, vol.~76, Princeton University Press, Princeton, NJ,
  1974.

\bibitem{Ndi07}
C.~B. Ndiaye, \emph{Constant ${Q}$-curvature metrics in arbitrary dimension},
  J. Funct. Anal. \textbf{251} (2007), no.~1, 1--58.

\bibitem{Pag78}
D.~N. Page, \emph{A compact rotating gravitational instanton}, Phys. Lett. B
  \textbf{79} (1978), no.~3, 235--238.

\bibitem{PagPop87}
D.~N. Page and C.~N. Pope, \emph{Inhomogeneous {E}instein metrics on complex
  line bundles}, Class. Quantum Grav. \textbf{4} (1987), 213--225.

\bibitem{Ste99}
N.~Steenrod, \emph{The topology of fibre bundles}, {R}eprint of the 1957 ed.,
  Princeton Landmarks in Mathematics, Princeton University Press, Princeton,
  NJ, 1999.

\bibitem{WanWan98}
J.~Wang and M.~Y. Wang, \emph{Einstein metrics on ${S}^2$-bundles}, Math. Ann.
  \textbf{310} (1998), no.~3, 497--526.

\bibitem{Wan98}
M.~Y. Wang, \emph{Einstein metrics from symmetry and bundle constructions},
  Surveys in differential geometry: essays on Einstein manifolds, Int. Press,
  Boston, MA, 1999, pp.~287--325.

\bibitem{WanZil90}
M.~Y. Wang and W.~Ziller, \emph{Einstein metrics on principal torus bundles},
  J. Differential Geom. \textbf{31} (1990), no.~1, 215--248.

\end{thebibliography}

\end{document}